\documentclass[graybox,envcountsect,envcountsame]{svmult}

\usepackage{mathptmx}       
\usepackage{helvet}         
\usepackage{courier}        
\usepackage{type1cm}        
\usepackage{makeidx}         
\usepackage{graphicx}        
\usepackage{multicol}        
\usepackage[bottom]{footmisc}

\makeindex             

\usepackage{amssymb,amsmath,mathrsfs,graphics,stmaryrd,enumitem,tikz,bm,chngcntr}

\newcommand{\defterm}[1]{\textbf{#1}}

\DeclareMathOperator{\covol}{covol}
\DeclareMathOperator{\colspace}{colspace}
\DeclareMathOperator{\im}{im}
\DeclareMathOperator{\Id}{Id}
\DeclareMathOperator{\nullity}{nullity}

\DeclareMathOperator{\rank}{rank}

\DeclareMathOperator{\Star}{star}
\DeclareMathOperator{\Sig}{Sig}
\DeclareMathOperator{\Cut}{Cut}
\DeclareMathOperator{\Flow}{Flow}
\DeclareMathOperator{\Tor}{Tor}

\newcommand{\0}{\emptyset}
\newcommand{\bd}{\partial}
\newcommand{\cbd}{\partial^*}
\newcommand{\dju}{\mathbin{\dot{\cup}}}

\newcommand{\eqdef}{\overset{{\rm def}}{=}}
\newcommand{\into}{\hookrightarrow}
\newcommand{\isom}{\cong}
\newcommand{\padup}{\rule{0mm}{5mm}}
\newcommand{\pdet}{\boldsymbol{\lambda}}
\newcommand{\sm}{\setminus}
\newcommand{\st}{~|~}
\newcommand{\tor}{{t}}
\newcommand{\Tree}{\mathscr{T}}
\newcommand{\x}{\times}

\newcommand{\cc}{\mathbf c}

\newcommand{\vv}{\mathbf{v}}
\newcommand{\ww}{\mathbf{w}}
\newcommand{\xx}{\mathbf{x}}
\newcommand{\yy}{\mathbf{y}}

\newcommand{\C}{\mathcal{C}}
\newcommand{\F}{\mathcal{F}}

\newcommand{\II}{\mathcal{I}}
\newcommand{\LL}{\mathcal{L}}

\newcommand{\OO}{\mathcal{O}}

\newcommand{\Qq}{\mathbb{Q}}
\newcommand{\Rr}{\mathbb{R}}
\newcommand{\Ss}{\mathbb{S}}
\newcommand{\Zz}{\mathbb{Z}}

\newcommand{\HH}{H}
\newcommand{\betti}{\beta}

\newcommand{\Ltot}{L^{\rm tot}} 
\newcommand{\Lud}{L^{\rm ud}} 
\newcommand{\Ldu}{L^{\rm du}} 
\newcommand{\Lalg}{L^{\rm alg}} 

\begin{document}

\title*{Simplicial and Cellular Trees}

\author{Art M.\ Duval \and Caroline J.\ Klivans \and Jeremy L. Martin}
\authorrunning{A.M.\ Duval, C.J.\ Klivans and J.L.\ Martin}
\institute{Art M.\ Duval \at Department of Mathematical Sciences, University of Texas at El Paso, \email{aduval@utep.edu}
\and Caroline J.\ Klivans \at Division of Applied Mathematics and Department of Computer Science, Brown University, \email{klivans@brown.edu}
\and Jeremy L. Martin \at Department of Mathematics, University of Kansas, \email{jlmartin@ku.edu}}
\maketitle

\abstract*{Much information about a graph can be obtained by studying 
its spanning trees.  On the other hand, a graph can be regarded as a
1-dimensional cell complex, raising the question of developing a
theory of trees in higher dimension.  As observed first by Bolker,
Kalai and Adin, and more recently by numerous authors, the fundamental
topological properties of a tree --- namely acyclicity and
connectedness --- can be generalized to arbitrary dimension as the
vanishing of certain cellular homology groups.  This point of view is
consistent with the matroid-theoretic approach to graphs, and yields
higher-dimensional analogues of classical enumerative results
including Cayley's formula and the matrix-tree theorem.  A subtlety of
the higher-dimensional case is that enumeration must account for the
possibility of torsion homology in trees, which is always trivial for
graphs.  Cellular trees are the starting point for further
high-dimensional extensions of concepts from algebraic graph theory
including the critical group, cut and flow spaces, and discrete
dynamical systems such as the abelian sandpile model.}

\section{Introduction} \label{section:Introduction}

How can redundancy be eliminated from a network, and what can be said about the resulting substructures?  In graph-theoretic terms, these substructures are spanning trees.  The collection of all spanning trees of a graph has the structure of a matroid basis system; this observation connects trees to algebraic combinatorics and explains why many graph algorithms can be made computationally efficient.  
The number of spanning trees of a graph measures its complexity as a network, and there are classical and efficient linear-algebraic tools for calculating this number, as well as for enumerating trees more finely.  The algebra of trees is key in studying dynamical systems on graphs, notably the abelian sandpile model (known in other forms as the chip-firing game or dollar game), whose possible states are encoded by a group of size equal to the complexity of the underlying graph.

This chapter is about the more recent theory of spanning trees in \emph{cell complexes}, which are natural higher-dimensional analogues of graphs. The theory for cell complexes parallels the graph-theoretical version in many ways, including the connection to matroid theory.  However, higher-dimensional spaces can have much richer topology, which complicates the algebraic and enumerative parts of the story in a very concrete way.  It turns out that the simple number of spanning trees of a cell complex is not a well-behaved invariant, and it is better to account for topological complexity in the enumeration.  On the other hand, many essential algebraic tools for working with spanning trees of a graph do extend well to arbitrary dimension.

Sections~\ref{section:Introduction} and~\ref{section:trees-and-forests} give an overview of the subject, including the history of the subject of higher-dimensional trees and emphasizing the geometry and topology of trees.  As far as possible, we have attempted to make these sections self-contained and accessible to readers not familiar with the machinery of algebraic topology, although a certain amount of technical material is unavoidable.  The next two sections are more technical: Section~\ref{section:enumeration} describes the various extensions of the matrix-tree theorem from graphs to cell complexes, and Section~\ref{section:specific-enumeration} examines individual classes of cell complexes whose tree enumerators are understood.  Finally, Section~\ref{section:open-problems} covers new directions, including extensions of algebraic graph theory and connections to critical groups, cuts and flows, simplicial decompositions, and matroids.

\subsection{Higher-dimensional trees} \label{subsection:intro-high-dim-trees}
 The higher-dimensional trees considered here were first studied by Bolker \cite{Bolker}, in the context of transportation polytopes, which are solution sets to certain combinatorial optimization problems.  The vertices of transportation polytopes are spanning trees of \emph{complete colorful complexes}, which are higher-dimensional generalizations of complete bipartite graphs.  Bolker pinpointed the basic homological conditions on spanning trees (essentially Definition~\ref{define-cellular-tree} below).
Kalai~\cite{Kalai} considered spanning trees of \emph{complete simplicial
complexes}, which generalize complete graphs, and obtained a beautiful enumerative formula generalizing Cayley's formula. Soon after, Adin~\cite{Adin} gave a corresponding formula for counting spanning trees in complete colorful complexes.
More recently, many authors, including Duval, Klivans and Martin~\cite{DKM1,DKM2,DKM4}, Lyons~\cite{Lyons} and Petersson~\cite{Petersson}, independently observed that Bolker
and Kalai's definitions could be extended to general cell
complexes to develop a broad theory of higher-dimensional trees.

Algebraically speaking, one can define a spanning tree of a suitably connected cell complex to be any subcomplex generated by faces corresponding to a linearly independent set of columns of a certain matrix (see Definition~\ref{define-cellular-forest}).
However, this definition does not really reflect the topology of higher-dimensional trees.  By analogy, a spanning tree $T$ of a connected graph $G$ can be defined as a collection of edges which correspond to a column basis of the oriented incidence matrix of $G$.   But this definition does not directly capture what a tree looks like, and the most familiar definitions are more topological and combinatorial in nature:
\begin{itemize}[]
\item $T$ is connected and acyclic;
\item $T$ is connected and $|T|=n-1$ (where $n$ is the number of vertices);
\item $T$ is acyclic and $|T|=n-1$;
\item $T$ is a maximal acyclic subgraph of $G$;
\item $T$ is a minimal connected subgraph of $G$.
\end{itemize}
Similarly, for a $d$-dimensional cell complex $X$ with codimension-1 Betti number $\betti_{d-1}(X)=0$ (the analogue of connectedness for graphs), a subcomplex $T\subseteq X$ is a spanning tree if it satisfies any of the following equivalent properties:
\begin{itemize}
\item $\betti_d(T)=\betti_{d-1}(T)=0$;
\item $\betti_d(T)=0$ and $|T_d|=|X_d|-\betti_d(X)$;
\item $\betti_{d-1}(T)=0$ and $|T_d|=|X_d|-\betti_d(X)$;
\item $T$ is maximal among the spanning subcomplexes of $X$ with $\betti_d(T)=0$;
\item $T$ is minimal among the spanning subcomplexes of $X$ with $\betti_{d-1}(T)=0$.
\end{itemize}
A more complete list of characterizations, including the more general notion of spanning forests when $\betti_{d-1}(X)>0$, is given in Proposition~\ref{characterize-forest}.  Definitions and notation for the terms above are also given in Section~\ref{section:trees-and-forests}.

Many cellular spanning trees ``look like'' trees.  For instance, a $2$-dimensional simplicial complex whose dual graph forms a tree automatically satisfies the conditions above; see the first complex in Figure~\ref{fig:trees}.  On the other hand, many higher-dimensional trees look less familiar.  For instance, if $X$ is any cellular sphere
(such as the boundary of a polytope), then its spanning trees are the contractible subcomplexes formed by removing a single maximal face; this is the higher-dimensional generalization of the fact that the spanning trees of a cycle graph are the paths formed from it by deleting one edge.  Of particular interest (and the source of much
difficulty) is the existence of trees with torsion homology, which can occur in any dimension $\geq 2$.  Already in dimension~$2$, trees can have finite non-trivial homology.  The smallest simplicial example is the standard triangulation of the real projective plane $\Rr P^2$, which arises as a spanning tree of the $2$-dimensional skeleton of the $6$-vertex simplex.  See Examples~\ref{RP2} and~\ref{K62}. 

\subsection{Tree numbers} \label{subsection:intro-enumeration}

The subject of enumerating spanning trees traditionally dates back to Kirchhoff's 1847 work on electrical circuits~\cite{Kirchhoff}.   Excellent general references with more historical details include Moon~\cite{Moon} and Stanley~\cite[pp.~54--69]{EC2}. 
In modern terms, let~$G$ be a connected graph on~$n$ vertices whose Laplacian matrix~$L=L(G)$ has nonzero eigenvalues $\lambda_1,\dots,\lambda_{n-1}$, and let $\Tree(G)$ denote its set of spanning trees.  Then
the \emph{matrix-tree theorem} (or \emph{Kirchhoff's theorem}) states that
\begin{equation} \label{MTT-orig}
\tau(G) = \det\tilde L = \frac{\lambda_1\cdots\lambda_{n-1}}{n}
\end{equation}
where $\tau(G)=|\Tree(G)|$ and $\tilde L$ denotes the reduced Laplacian obtained from~$L$ by deleting the~$i^{th}$ row and column (for any $i$).
Thus tree enumeration is closely tied in with the broad subject of spectral graph theory; see, e.g., \cite{SpectraBook, Chung-book}.
A most important special case is the formula
\begin{equation} \label{cayley}
\tau(K_n)=n^{n-2}
\end{equation}
for the number of spanning trees of the complete graph~$K_n$, known as  \emph{Cayley's formula}, although it was in fact first observed by Sylvester~\cite{Sylvester} and first proved, in an equivalent form, by Borchardt~\cite{Borchardt}.
Cayley's formula can be deduced from the matrix-tree theorem by finding an explicit eigenbasis for $L(K_n)$.  Independent combinatorial proofs
include the well-known Pr\"ufer code \cite{Prufer} and Joyal's bijection using permutations~\cite[pp.~15-16]{Joyal}; see also~\cite[pp.~25--28]{EC2}. 

There are other families of graphs with nice spanning tree counts which can be obtained either algebraically (using the matrix-tree theorem) or combinatorially (via a bijection).  For example, the number of spanning trees of  the complete bipartite graph $K_{m,n}$ is
\begin{equation} \label{count-Kmn}
\tau(K_{m,n})=n^{m-1}m^{n-1}
\end{equation}
as discovered by several authors \cite{Austin,Simmonard,Fiedler-Sedlacek,Scoins}, and arises in the {transportation problem} of combinatorial optimization \cite{Klee-Witzgall}.

Moving to higher dimensions, in Bolker's original work, he
observed that it is ``tempting on numerological grounds'' to generalize the well-known formulas~\eqref{cayley} and~\eqref{count-Kmn}, but also was aware that the obvious generalizations were incorrect because of \defterm{torsion trees}, i.e., $d$-dimensional
trees $T$ with $\HH_{d-1}(T;\Zz)$ finite but nontrivial.
Kalai~\cite{Kalai} was the first to pinpoint the role played by torsion.  He showed that the correct analogue of the tree count for a $d$-dimensional simplicial complex $X$ is not simply the cardinality of the set $\Tree(X)$ of simplicial spanning trees, but rather the quantity
\begin{equation} \label{define-tau}
\tau(X)=\sum_{T\in\Tree(X)}|\HH_{d-1}(T)|^2
\end{equation}
where $\HH_{d-1}(T)$ denotes the $(d-1)$-dimensional reduced homology group of $T$ over $\Zz$.  The summands in~\eqref{define-tau}
arise naturally in the expansion of the combinatorial Laplacian matrix of~$X$ using the Binet-Cauchy formula.  Note that torsion has no effect on spanning tree enumeration in graphs; every topological space has torsion-free $0^{\rm th}$ homology group, so when $X$ is a graph the right-hand side of~\eqref{define-tau} reduces to the number of spanning trees.  (On the other hand, torsion does arise in enumerative problems equivalent or related to tree enumeration; see, e.g.,~\cite{JLM-Cyclotomic,Beck-Hosten,BBGM, Moci1}.)
Additional evidence that $\tau(X)$ is the right way to count trees is Kalai's striking generalization of Cayley's formula:
\begin{equation} \label{kalai}
\tau(K_{n,d})=n^{\binom{n-2}{d}}
\end{equation}
where $K_{n,d}$ denotes the $d$-dimensional skeleton of an $n$-vertex simplex.  Subsequently, Adin~\cite{Adin} used the same ideas to give a formula for the torsion-weighted count of a complete colorful complex, generalizing~\eqref{count-Kmn}.  Both formulas were precisely those predicted by Bolker.  
Kalai also gave a weighted version of~\eqref{kalai} enumerating spanning trees by their degree sequences~\cite[Thm.~$3'$]{Kalai} and showed that torsion accounts for most of the quantity $\tau(X)$: as Kalai put it, spanning trees of $K_{n,d}$ are ``on the average, far from being [$\Zz$-]acyclic''~\cite[\S5]{Kalai}.

The torsion-weighted spanning tree count of a cell complex $X^d$ can be computed from its combinatorial Laplacian operator $L$.  This essential fact, the \emph{cellular matrix-tree theorem}, has appeared in various forms in~\cite{DKM1, DKM2, Lyons, Petersson,  CCK}, all building on the ideas originating in~\cite{Kalai}
and~\cite{Adin}.  We will describe the different guises of this theorem, and their interconnections, in Section~\ref{section:enumeration} below.  As in the
graphical case, $\tau(X)$ may be expressed via a determinant or via
eigenvalues.  Specifically, for a $d$-dimensional complex $X$, the torsion-weighted tree count of $d$-dimensional trees is given by
\[\tau_d(X) \, = \, c \cdot \det \tilde L \, = \, c'\cdot  \frac{\lambda_1 \cdots \lambda_{k}}{\tau_{d-1}(X)}\]
where $c$ and $c'$ are certain correction factors arising from torsion homology; $\tilde L$ is the reduced Laplacian formed from~$L$ by removing the rows and columns corresponding to a codimension-1 spanning tree; and the
$\lambda_i$ are the non-zero eigenvalues of~$L$.

In practice, for many families of complexes of combinatorial interest (e.g., Cohen-Macaulay complexes), the homological
correction factors $c$ and $c'$ disappear.  In this case, the weighted tree count is given precisely by a determinant, as in Kalai's enumeration of skeletons of simplices.  In general, the
eigenvalue formulation of the matrix-tree theorem is recursive; the
count for $d$-dimensional trees is expressed in terms of the
count for $(d-1)$-dimensional trees. In the case that $c'=1$, the recursion yields the alternating-product expression
\[\tau_d(X) = \frac{ \pdet_1 \cdot \pdet_3 \cdot \pdet_5 \cdots}{\pdet_0 \cdot \pdet_2 \cdot \pdet_4 \cdots}\]
where $\pdet_i$ is the product of the non-zero eigenvalues of the $i$th Laplacian.

Laplacians of simplicial complexes first appeared in the work of Eckmann~\cite{Eckmann} and Dodziuk
and Patodi~\cite{Dodziuk}.  They were studied as combinatorial objects
by Friedman~\cite{Friedman} and Friedman and Hanlon~\cite{Friedman-Hanlon}, who in particular considered the phenomenon of
integral Laplacian eigenvalues.  Subsequently, many combinatorially
significant classes of complexes have been shown to have all integer
Laplacian eigenvalues, such as chessboard
complexes~\cite{Friedman-Hanlon}, matroid complexes~\cite{KRS},
shifted complexes~\cite{Duval-Reiner}, and hypercubes~\cite{DKM2}.
The strong connection between Laplacian eigenvalues and tree
enumeration suggests that these families might support nice spanning
tree counts with combinatorial interpretations.  Indeed, there are many successes in this direction including enumerators for shifted complexes in terms of critical pairs~\cite{DKM1} and for matroid complexes in terms of flats and the $\beta$-invariant~\cite{KookLee-FSTNMC}.  We will look at many of these families individually in Section~\ref{section:specific-enumeration}.

\subsection{Beyond basic definitions and enumeration} \label{subsection:intro-beyond} 

The theory of trees has recently emerged as a natural aspect of the combinatorics of cellular spaces.  We describe several new directions in Section~\ref{section:open-problems}.  Many fundamental objects of algebraic graph theory have been generalized to higher dimension, including critical groups and generalized sandpile/chip-firing models (\S\ref{subsection:critical}) and cuts and flows (\S\ref{subsection:cutflow}).  Simplicial trees provide an inroad to studying theorems and conjectures about decompositions of simplicial complexes (\S\ref{subsection:decomp}).  Finally, the study of cellular matroids has led to new topological invariants on cellular manifolds and applications to probability (\S\ref{subsection:trees-and-matroids}).
\section{Trees and forests: from graphs to cell complexes} \label{section:trees-and-forests}

The definitions of higher-dimensional trees and forests may be approached from various directions, including linear algebra and matroid theory, but ultimately some topological machinery is necessary to appreciate the ideas fully.  
We will attempt to give an informal, self-contained description of those concepts that we need such as cell complexes and homology, that will be accessible to a reader  willing to endure a little black-box algebra.  For the precise details, the reader should consult a standard reference such as Hatcher~\cite{Hatcher}.

\subsection{Trees and forests in graphs} \label{subsection:tree-forest-graph}

Let $G=(V,E)$ be a finite graph, with $|V|=n$ and $|E|=m$.  For each edge $e$ on vertices $v$ and $w$, fix an orientation $e=v\to w$.  The \defterm{signed incidence matrix} $\bd_G$ is then defined as the $n\x m$ matrix with entries
\[\bd_{v,e}=\begin{cases}
1 & \text{ if $e=u\to v$ for some $u\neq v$},\\
-1 & \text{ if $e=v\to w$ for some $w\neq v$},\\
0 & \text{ otherwise.} \end{cases}\]
In particular, the $e^{th}$ column of $\bd$ is zero if and only if $e$ is a loop.
We now recall several well-known equivalent definitions of trees and forests from graph theory.

\begin{definition} \label{defn:graph-tree-forest}
The graph $G$ is a \defterm{tree} if it satisfies any of the following equivalent properties:
\begin{enumerate}[label=(\alph*)]
\item $G$ is connected and acyclic.
\item $G$ is connected and $m=n-1$.
\item $G$ is acyclic and $m=n-1$.
\item\label{unique-path-trees} Every pair of vertices in $G$ is connected by exactly one path.
\item $G$ is a maximal acyclic graph.
\item $G$ is a minimal connected graph.
\item For any (hence every) orientation of $G$, the columns of $\bd_G$ are a basis for the space 
\(\Rr^n_0=\{v\in\Rr^n \st v_1+\cdots+v_n=0\}\).
\end{enumerate}
Meanwhile, $G$ is a \defterm{forest} if it satisfies any of the following equivalent properties:
\begin{enumerate}[resume,label=(\alph*)]
\item $G$ is acyclic.
\item Every pair of vertices in $G$ is connected by at most one path.
\item The columns of the incidence matrix $\bd_G$ of $G$ are linearly independent.
\item\label{forest-is-trees} Every connected component of $G$ is a tree.
\end{enumerate}
\end{definition}

In order to extend the ideas of trees and forests, we first have to describe the objects that will serve as higher-dimensional analogues of graphs.

\subsection{Cell complexes} \label{subsection:cell-complexes}

A \defterm{cell complex} $X$ is a topological space that is the disjoint union of \emph{cells}, each homeomorphic to the interior of a $k$-dimensional disk
for some~$k$, attached together in a locally reasonable manner.  In particular, the closure of any cell is a union of cells.  The set of $k$-cells is denoted $X_k$, and the \defterm{$k$-skeleton} $X_{\leq k}$ is the union of all cells of dimension $\leq k$.  A cell not contained in the closure of any other cell is called a \defterm{facet}, and
a complex whose facets all have the same dimension is called \defterm{pure}.  Every subset $S$ of cells generates a subcomplex $X_S=\overline{\bigcup_{\sigma\in S}\sigma}$.  (We often blur the distinction between a set of cells and the subcomplex it generates.)
The dimension of $X$ is the largest dimension of one of its cells; the notation $X^d$ indicates a cell complex of dimension~$d$.  In particular, a graph is exactly a cell complex of dimension 1.

The combinatorial data about how cells are attached is recorded by a sequence of linear maps $\bd_k=\bd_k(X):\Rr^{X_k}\to\Rr^{X_{k-1}}$ for $0\leq k\leq d$.  The map $\bd_1$ is just the (signed) incidence matrix of the 1-skeleton graph, and the other maps are higher-dimensional analogues of the signed incidence matrix of a graph.  For $k>0$, each $\bd_k$ may be represented as a matrix with rows indexed by $(k-1)$-cells $\rho$ and columns indexed by $k$-cells $\sigma$; the entry $\bd_{\rho,\sigma}$ is an integer that records, roughly speaking, the number of times that $\sigma$ is wrapped around $\rho$, counted with orientation.  If $\rho$ is not contained in the closure of $\sigma$, then $\bd_{\rho,\sigma}=0$.  In particular, if some $(k-1)$-cell $\rho$ is not in the closure of any $k$-cell, then the corresponding row of $\bd$ is zero.  We formally define $\bd_0:\Rr^{X_0}\to\Rr$ to be the map that sends every point to $1\in\Rr$.

It is a standard fact that $\im\bd_{k+1}\subseteq\ker\bd_k$ for all $k$, and the quotient
$\HH_k(X;\Rr)=\ker\bd_k/\im\bd_{k+1}$
is a topological invariant of $X$ called its \defterm{$k^{\rm th}$ (reduced) homology group}.  The \defterm{$k^{\rm th}$ (reduced) Betti number} is
$\betti_k(X)=\dim_\Rr\HH_k(X;\Rr)$. (All homology groups and Betti numbers appearing here are reduced and we will use the notation $\HH_k$ and $\betti_k$ in place of the usual $\tilde{H}_k$ and $\tilde{\beta}_k$ for ease of readability.)  The $0^{\rm th}$ reduced Betti number is one less than the number of connected components, so in particular $\betti_0(X)=0$ if and only if $X$ is connected.  If $G=(V,E)$ is a graph with $c$ connected components, then $\betti_1(G)=|E|-(|V|-c)$ is the size of the complement of a maximal spanning forest (sometimes called the \emph{cyclomatic number}).

\defterm{Simplicial complexes} are the cell complexes most familiar to combinatorialists.  An (abstract) simplicial complex is a set $\Delta$ of subsets of $[n]=\{1,2,\dots,n\}$, such that whenever $\sigma\in\Delta$ and $\rho\subseteq\sigma$, then $\rho\in\Delta$.  Thus a simplicial complex is completely determined by its set of facets.
An abstract simplicial complex can be regarded as a topological space by associating $\sigma=\{i_1,\dots,i_k\}\in\Delta$ with the interior of the convex hull of $\{e_{i_1},\dots,e_{i_k}\}$ in $\Rr^n$, where $e_j$ is the $j^{th}$ standard basis vector.  For readers familiar with hypergraphs, a simplicial complex is just a hypergraph in which every subset of a hyperedge is also a hyperedge.

Simplicial complexes are \defterm{regular} cell complexes, which means that every attaching map is locally a homeomorphism; i.e., one-to-one.  In particular, for a regular cell complex, every entry of every matrix $\bd_k$ is in $\{0, \pm 1\}$.

\begin{example} \label{example:G}
A graph $G$ is a 1-dimensional cell complex whose attaching map $\bd_1$ is just the signed vertex-edge incidence matrix.  As a cell complex, $G$ is regular if $G$ has no loops, and simplicial if in addition $G$ has no parallel edges.
\end{example}

\begin{example} \label{example:B}
Gluing two hollow tetrahedra along a common triangle 123 produces a cell complex, the \defterm{equatorial bipyramid}, with seven 2-cells (triangles), nine 1-cells (edges), and five 0-cells (vertices).  This is an example of a simplicial complex, with facets
$\{123,124,125,124,134,234,235\}$
where 1,2,3 are the ``equatorial'' vertices and 4,5 are the ``poles''.  The collection of facets would often be considered a $3$-uniform hypergraph; i.e., a hypergraph where all edges have size $3$.
\end{example}

\begin{example} \label{example:RP2}
The real projective plane $\Rr P^2$ can be realized as a cell complex with one cell $\sigma_k$ of each dimension $k$ for $k\in\{0,1,2\}$.  The 1-skeleton $\sigma_1\cup \sigma_0$ forms the line at infinity; the attaching map $\bd_1$ is the matrix $[0]$ because $\sigma_1$ is a loop.  Meanwhile, the attaching map $\bd_2$ is the matrix $[2]$, representing the fact that $\sigma_2$ is wrapped twice around $\sigma_1$. The smallest triangulation of $\Rr P^2$ (i.e., the smallest simplicial complex homeomorphic to $\Rr P^2$) can be obtained by identifying antipodal faces on an icosahedron.  See Figure~\ref{fig:RP2}.  
\end{example}

\subsection{Trees and forests in cell complexes} \label{subsection:trees-forests-cell-complexes}
 
We can now define the cellular analogues of trees and forests.

\begin{definition} \label{define-cellular-forest}
A $d$-dimensional cell complex $X$  is a \defterm{$d$-forest} if the columns of $\bd_d(X)$ are linearly independent.
A \defterm{spanning $d$-forest} of $X$ is a $d$-forest $F\subseteq X$ such that $F_{\leq d-1}=X_{\leq d-1}$.
\end{definition}

Here the term ``spanning'' refers to the requirement that the $(d-1)$-skeleton of $F$ equals that of $X$, just as a spanning subgraph of a graph $G$ is a subgraph that includes all vertices.
From a homological perspective, an  equivalent condition for $X$ to be a $d$-forest is that $\HH_d(X;\Rr)=0$, because $\HH_d(X_{\leq d};\Rr)=\ker\bd_d(X)$.

\begin{definition} \label{define-cellular-tree}
A $d$-dimensional cell complex $X$  is a \defterm{$d$-tree} if the columns of $\bd_d(X)$ are linearly independent and $\rank\bd_d(X)=\nullity\bd_{d-1}(X)$.
A \defterm{spanning $d$-tree} of $X$ is a $d$-tree $T\subseteq X$ such that $T_{\leq d-1}=X_{\leq d-1}$.
\end{definition}
Equivalently, $X$ is a tree if $\HH_d(X;\Rr)=\HH_{d-1}(X;\Rr)=0$.

Note that every cell complex has at least one spanning forest (since its boundary matrix has at least one column basis), but a spanning forest of~$X$ is a spanning tree if and only if $\HH_{d-1}(X;\Rr)=0$.  In general, the conditions $\HH_d(X;\Rr)=0$ and $\HH_{d-1}(X;\Rr)=0$ should be regarded as the analogues of acyclicity and connectedness for
a $d$-dimensional cell complex.  Thus
Definitions~\ref{define-cellular-forest} and~\ref{define-cellular-tree} say that $d$-forests and $d$-trees are respectively ``$d$-acyclic'' and ``$d$-acyclic and $d$-connected'' cell complexes.

 What do cellular trees and forests look like?  A 1-dimensional forest is a graph that contains no cycles (subgraphs homeomorphic to a circle) and a $1$-dimensional tree contains no cycles and is connected (i.e., every pair of vertices can be seen as the endpoints of some $1$-dimensional path).  Correspondingly, a 2-dimensional cell complex~$X$ is a forest if~$X$ contains no subcomplex homeomorphic to a compact oriented 2-manifold.  Imagine a 2-dimensional cell complex as a network of bubbles: in order to make it into a forest,
one must remove enough 2-cells to pop all the bubbles, without removing any lower-dimensional cells.  Furthermore, $X$ is a 2-tree if every 1-cycle is the boundary of some disk consisting of 2-cells.  
Definitions~\ref{define-cellular-forest} and~\ref{define-cellular-tree} make sense even when $d=0$: a 0-dimensional forest is any set of vertices, and a 0-dimensional tree is just a single vertex, which is a 0-spanning tree of~$X$ if and only if $X_{\leq 1}$ is a connected graph.

\begin{figure}[ht]
\begin{center}
\includegraphics[width=4in]{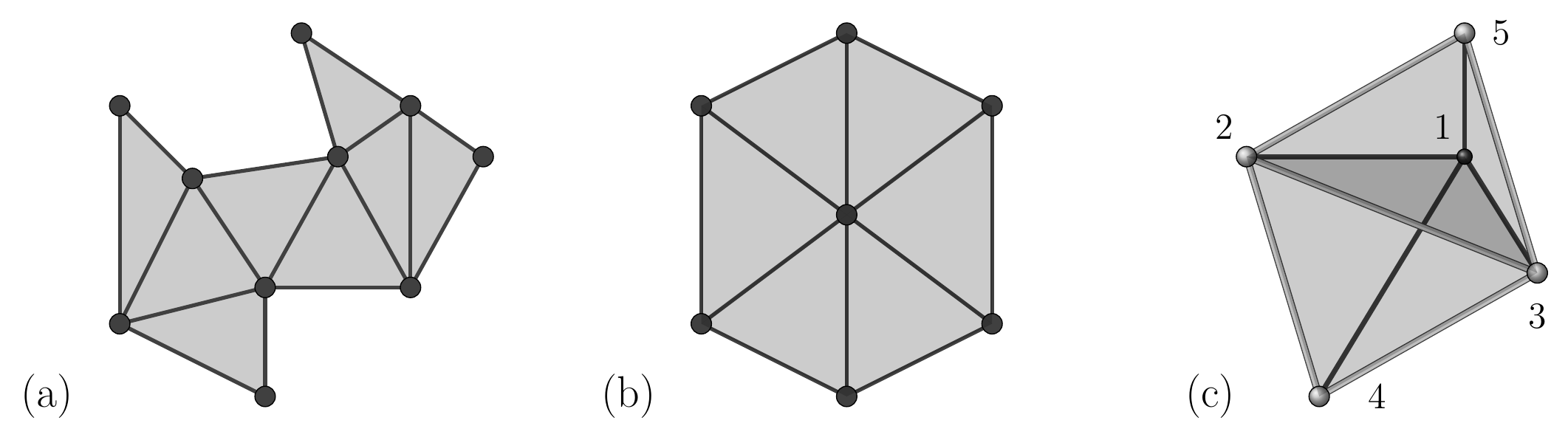}
\caption{Some two-dimensional trees.
\label{fig:trees}}
\end{center}
\end{figure}

\begin{example}\label{example:flat}
Figure~\ref{fig:trees} shows three 2-dimensional simplicial trees.  Complex~(a) is a triangulation of the disk whose dual graph is a tree.  Complex~(b) is also a triangulation of the
disk, but the dual graph forms a cycle, so it is a 2-tree with no ``leaf''.   Complex~(c) is a contractible simplicial complex (hence a tree) with facets 123, 124, 134, 125, 135.  It is a spanning tree of the equatorial bipyramid.
\end{example}

\begin{example} \label{manifold-trees}
If $X^d$ is homeomorphic to a compact oriented manifold (such as a $d$-sphere) then $\betti_d(X)=1$, and the kernel of $\bd_d$ is one-dimensional, generated by a linear relation that involves all its columns nontrivially.  Informally, $X$ has exactly one bubble, and it can be popped by removing any facet.
Deleting subsequent facets pokes extra holes in the surface.
Hence any spanning subcomplex obtained from~$X$ by removing one facet (resp., one or more facets) is a maximal spanning forest (resp., a spanning forest).
\end{example}

\begin{example} \label{bipyramid-tree}
Consider the equatorial bipyramid $B$ of Example~\ref{example:B}.  Topologically, it is a ``double bubble'', recorded homologically by its second Betti number $\betti_2(B)=2$.  In order to pop both bubbles, we can remove either (a) two triangles, one containing vertex 4 and one containing vertex 5, or (b) the middle triangle 123 plus any other triangle.  Thus $B$ has a total of $(3\x3)+6=15$ spanning trees, all of which are contractible simplicial complexes with 5 facets.  See Figure~\ref{fig:trees}(c) and Example~\ref{example:flat} for one example.
\end{example}

\begin{figure}[ht]
\begin{center}
\includegraphics[width=3in]{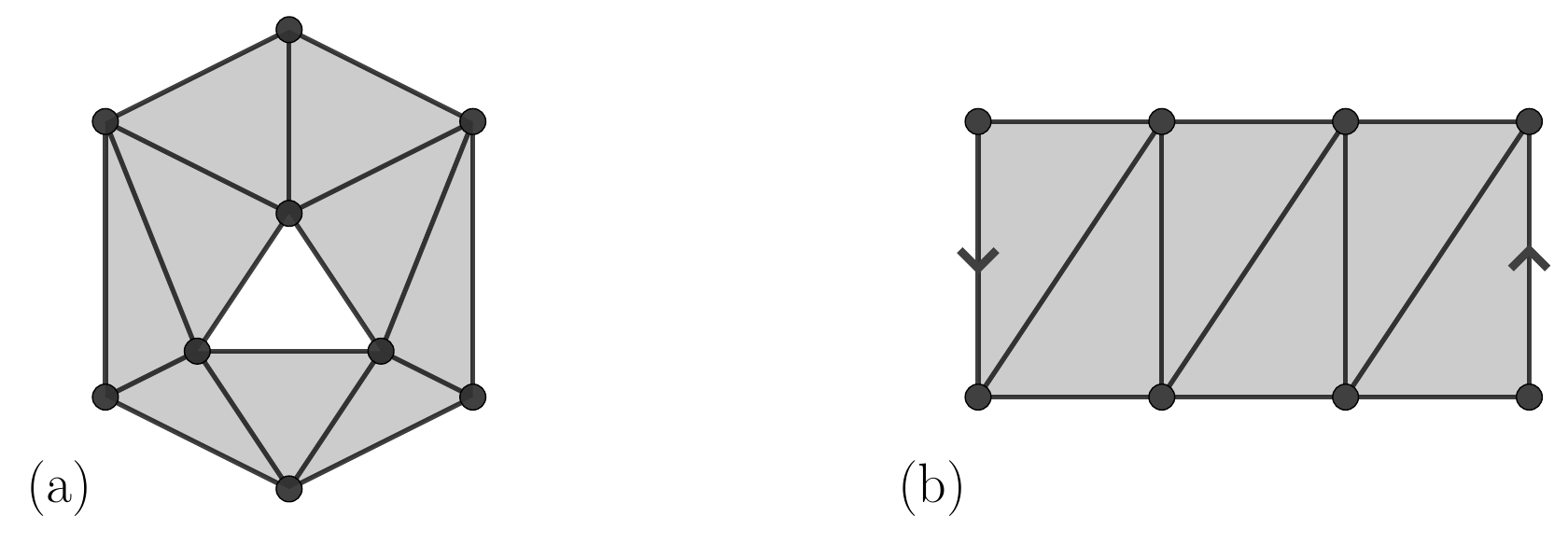}
\caption{Two-dimensional forests.
\label{fig:forests}}
\end{center}
\end{figure}

\begin{example}\label{example:forests}
Figure~\ref{fig:forests} shows two 2-dimensional simplicial forests, which triangulate the annulus and the M\"obius band.   In both cases $\betti_2=0$ but $\betti_1=1>0$ (since neither complex is simply connected), so each complex has itself as a spanning 2-forest, but has no spanning tree --- just like a non-connected acyclic graph.  On the other hand, unlike a 1-dimensional forest, neither complex is the disjoint union of 2-trees.
\end{example}

The connectedness condition  $\HH_{d-1}(X;\Rr)=0$ does not preclude the existence of \emph{torsion}.  Suppose that we regard the boundary matrix $\bd_k$ as a map on $\Zz$-modules $\bd_k^{\Zz}:\Zz^{X_k}\to\Zz^{X_{k-1}}$ (rather than on $\Rr$-vector spaces), and accordingly work with the integer homology groups $\HH_k(X;\Zz)=\ker\bd^\Zz_k/\im\bd^\Zz_{k+1}$.  While working over $\Zz$ does not change the ranks of the matrices, the condition $\rank\bd_d(X)=\nullity\bd_{d-1}(X)$ no longer implies that 
$\HH_{d-1}(X;\Zz)$ is the trivial group, merely that it is \emph{finite}.  This subtlety is concealed in the graph case $d=1$, since it is a fact that $\HH_0(X;\Zz)$ is free abelian for any space $X$, but is of central importance in higher dimension.

\begin{figure}[h]
\begin{center}\includegraphics[width=1in]{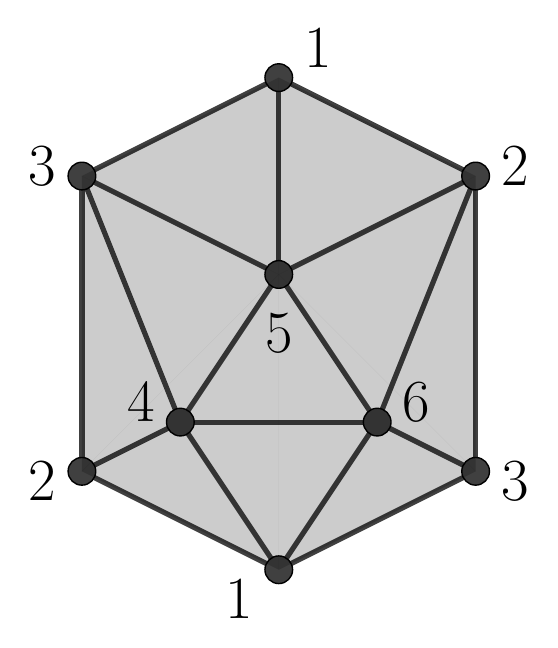}\end{center}
\caption{The triangulation of the real projective plane is a $2$-tree with non-trivial torsion.
\label{fig:RP2}}
\end{figure}

\begin{example} \label{RP2}  While every tree graph is a contractible topological space, this is not always true in higher dimension.  The simplest example is the real projective plane.  In the standard cell structure on $\Rr P^2$ described in Example~\ref{example:RP2}, $\im\bd_2=2\Rr=\ker\bd_1(X)$.
So $\HH_1(\Rr P^2;\Rr)=0$ and $\Rr P^2$ is a 2-tree (hence the only spanning 2-tree of itself).  On the other hand, $\HH_1(\Rr P^2;\Zz)\isom\Zz/2\Zz\neq0$, so $\Rr P^2$ is not contractible.  In the triangulation of $\Rr P^2$ shown in Figure~\ref{fig:RP2}, the cycle consisting of edges 12, 13, and 23 is a representative of the non-trivial torsion homology class.
\end{example}

\begin{example}\label{K62}
Let $K_{6,2}$ be the 2-dimensional skeleton of the 6-vertex simplex, i.e., the
2-dimensional simplicial complex on vertex set $V=[6]$ whose facets are the twenty 3-element subsets of $V$.  Let $T_1$ be the spanning subcomplex of~$K_{6,2}$ whose facets are the ten 3-element sets containing~$v_1$.  Then $T_1$ is contractible (it is a cone), hence is a simplicial spanning tree of~$K_{6,2}$ with $\HH_1(T_1;\Zz)=0$.  On the other hand, $\Rr P^2$ can be triangulated with six vertices (as in Figure~\ref{fig:RP2}), and any such triangulation $T_2$ is also a spanning tree of~$K_{6,2}$, with $\HH_1(T_2;\Zz)\isom\Zz/2\Zz$.  
\end{example}

Just as deleting an edge from a graph either decrements the cyclomatic number or increments the number of components, deleting facets of a simplicial complex $X^d$ either decrements $\betti_d$ or increments $\betti_{d-1}$.  Thus a  $d$-forest is a subcomplex $T\subset X$ obtained by removing enough facets so that $\betti_d(T)=0$, and a maximal  $d$-forest is one in which no facets have been removed unnecessarily, so that $\betti_{d-1}(T)=\betti_{d-1}(X)$.  If in addition $\betti_{d-1}(X)=0$, then $T$ is a spanning tree.  Thus the vanishing of $\betti_d(T)$ and $\betti_{d-1}(T)$ are the higher-dimensional analogues of the graph-theoretic concepts of acyclicity and connectedness respectively.

The number of facets of every maximal  $d$-forest $T\subseteq X$ is $|T_d|=|X_d|-\betti_d(X)$, which is the analogue of the statement that a maximal spanning forest of a graph with $n$ vertices and $c$ components has $n-c$ edges. 

This discussion leads to the following characterizations of cellular spanning forests and trees, which generalize most of the conditions of Definition~\ref{defn:graph-tree-forest}.

\begin{proposition} \label{characterize-forest}
Let $X$ be a cell complex of dimension $d$ and let $T$ contain all faces of dimension $(d-1)$ or less of $X$.  Then $T$ is a maximal spanning forest of $X$ if it satisfies any of the following equivalent properties:
\begin{enumerate}[label=(\alph*)]
\item $\betti_d(T)=0$ and $\betti_{d-1}(T)=\betti_{d-1}(X)$.
\item $\betti_{d-1}(T)=\betti_{d-1}(X)$ and $|T_d|=|X_d|-\betti_d(X)$.
\item $\betti_d(T)=0$ and $|T_d|=|X_d|-\betti_d(X)$.
\item\label{exactlyonechain} Every element of $\ker\bd_{d-1}(X,T)$ is the boundary of exactly one $d$-chain in $T$.
\item $T$ is maximal among the spanning subcomplexes of $X$ with $\betti_d(T)=0$.
\item $T$ is minimal among the spanning subcomplexes of $X$ with $\betti_{d-1}(T)=\betti_{d-1}(X)$.
\item The columns of $\bd_d(T)$ are a vector space basis for $\colspace\bd_d(X)$.
\end{enumerate}
Meanwhile, $T$ is a spanning forest (not necessarily maximal) if
\begin{enumerate}[resume,label=(\alph*)]
\item $\betti_d(T)=0$.
\item\label{atmostonechain} Every element of $\ker\bd_{d-1}(X,T)$ is the boundary of at most one $d$-chain in $T$.
\item The columns of $\bd_d(T)$ are linearly independent.
\end{enumerate}
\end{proposition}
A maximal spanning forest is a spanning tree when $\betti_{d-1}(X)=0$.  In this case $\ker\bd_{d-1}(X,T)$ becomes $\ker\bd_{d-1}(X)$ in~\ref{exactlyonechain} and~\ref{atmostonechain}.

When $d=1$, these properties reduce to the familiar equivalent definitions of trees and forests given in Definition~\ref{defn:graph-tree-forest}, letter by letter.  One condition that does not carry over to higher dimension is~\ref{forest-is-trees} of Definition~\ref{defn:graph-tree-forest}, which defines a forest as a graph in which every connected component is a tree, it relies on the fact that the incidence matrix~$\bd_1(G)$ of a graph $G$ breaks up as a block sum of the incidence matrices of its components, of which there are $\betti_0(G)+1$.  On the other hand, there is no canonical way to decompose a $d$-dimensional cell complex~$X$, or the matrix~$\bd_d(X)$, into $\betti_{d-1}(X)+1$ pieces. Hence, $d$-forests are not simply disjoint unions of $d$-trees, as in the annulus and M\"obius band of Example~\ref{example:forests}.

\begin{remark} \label{ex:z-acyclic}
If $T\subseteq X$ is a spanning tree, then there is a surjection of finite groups $\HH_{d-1}(T;\Zz)\to\HH_{d-1}(X;\Zz)$.  In particular, if $X$ has a $\Zz$-acyclic spanning tree then $\HH_{d-1}(X;\Zz)=0$, but the converse does not hold.  For example, consider a cell complex with top boundary map $\begin{bmatrix}p&q\end{bmatrix}$ for any relatively prime integers $p,q>1$.
\end{remark}

\subsection{Rooted trees and forests} \label{subsection:rooted}

What is the higher-dimensional analogue of a \emph{rooted} tree or forest?  In dimension~1, a root of a tree graph is a vertex, while a root of a forest graph is a choice of a vertex from each connected component.  As before, in order to generalize these concepts to arbitrary dimension, we first state them in terms of linear algebra.  It is a standard fact that the nullvectors of a graph $G$ are just the characteristic vectors of its connected components; therefore, a root of a maximal spanning forest of $G$ can be regarded as the complement of a column basis of $\bd_G$, which motivates the following definition.

\begin{definition} \label{defn:root}
Let $X$ be a $d$-dimensional cell complex.  A \defterm{rooted spanning forest} of $X$ is a pair $(F,R)$, where
$F$ is a spanning forest of $X$; $R$ is a $(d-1)$-dimensional subcomplex of $X$ such that $R_{\leq d-2} = X_{\leq d-2}$; and $\bd_{X\setminus R,F}$ is square and nonsingular, where $\bd_{X\setminus R,F}$ is
the submatrix of $\bd_d(X)$ with rows~$X \setminus R$ and columns~$F$.
\end{definition}

\begin{figure}[h]
\includegraphics[width=4.6in]{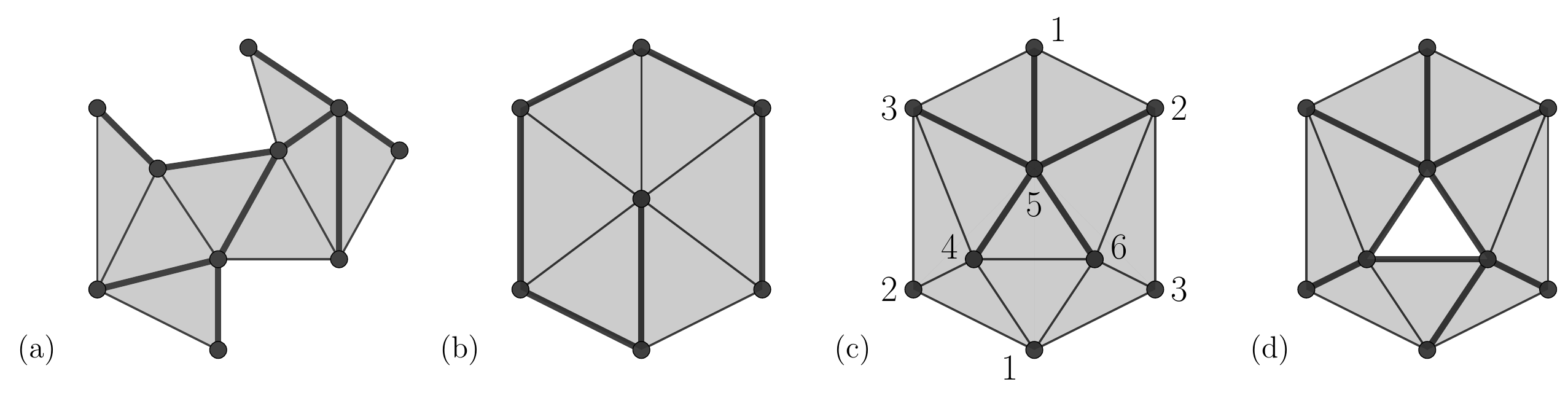}
\caption{Rooted trees and forests. \label{fig:rooted}}
\end{figure}

Figure~\ref{fig:rooted} shows rootings of some of the 2-trees encountered previously: the triangulations of the 2-disk (a,b), the projective plane (c), and the annulus (d).  The roots are the 1-dimensional subcomplexes indicated with thick lines.  In the first three cases $\betti_2(X)=0$, so a root of $X$ is simply a spanning tree of the
codimension-$1$ skeleton.  The annulus has $\betti_1=1$, and its root is a unicyclic graph rather than a spanning 1-tree.

In the case that $(F,R)$ is a maximal rooted spanning forest, the root
$R$ was called a ``relatively acyclic complex'' in \cite{DKM4},
because the inclusion map $R\into X$ induces
isomorphisms $\HH_i(R;\Rr)\to\HH_i(X;\Rr)$ for all $i<d$, a condition
equivalent  to the vanishing of
the relative homology groups $\HH_i(X,R;\Rr)$ for all $i<d$.
Lyons~\cite{Lyons} used the term ``cobase'' for the complement of a
root, referring to the fact that a $i$-cobase is a row basis for the
coboundary map $\bd_{i+1}^*$.

\begin{definition} \label{defn:orientation}
Let $X$ be a $d$-dimensional cell complex.  A \defterm {directed rooted spanning forest} of $X$ is a triple $(F,R,\OO)$ where $(F,R)$ is a rooted spanning forest of $X$, and $\OO: X_{d-1} \sm R_{d-1} \rightarrow F_{d}$ is a bijection
such that $\rho \subset \OO(\rho)$ for all $f$.  That is, $\OO$ pairs
each codimension-$1$ face $\rho\not\in R$ with a facet that contains it.
\end{definition}

When $d=1$, this formalism regards an oriented
edge $u \to v$ not as an ordered pair of vertices, but as the pairing of vertex~$u$ with edge~$uv$.  It is a standard fact that every rooting of a tree graph induces a unique orientation of its edges so that each edge is directed towards the root.
In higher dimensions, a root can admit more than one orientation; this subject is discussed more fully in \S\ref{subsection:decomp}.

\begin{figure}[h]
\begin{center}
\includegraphics[width=3in]{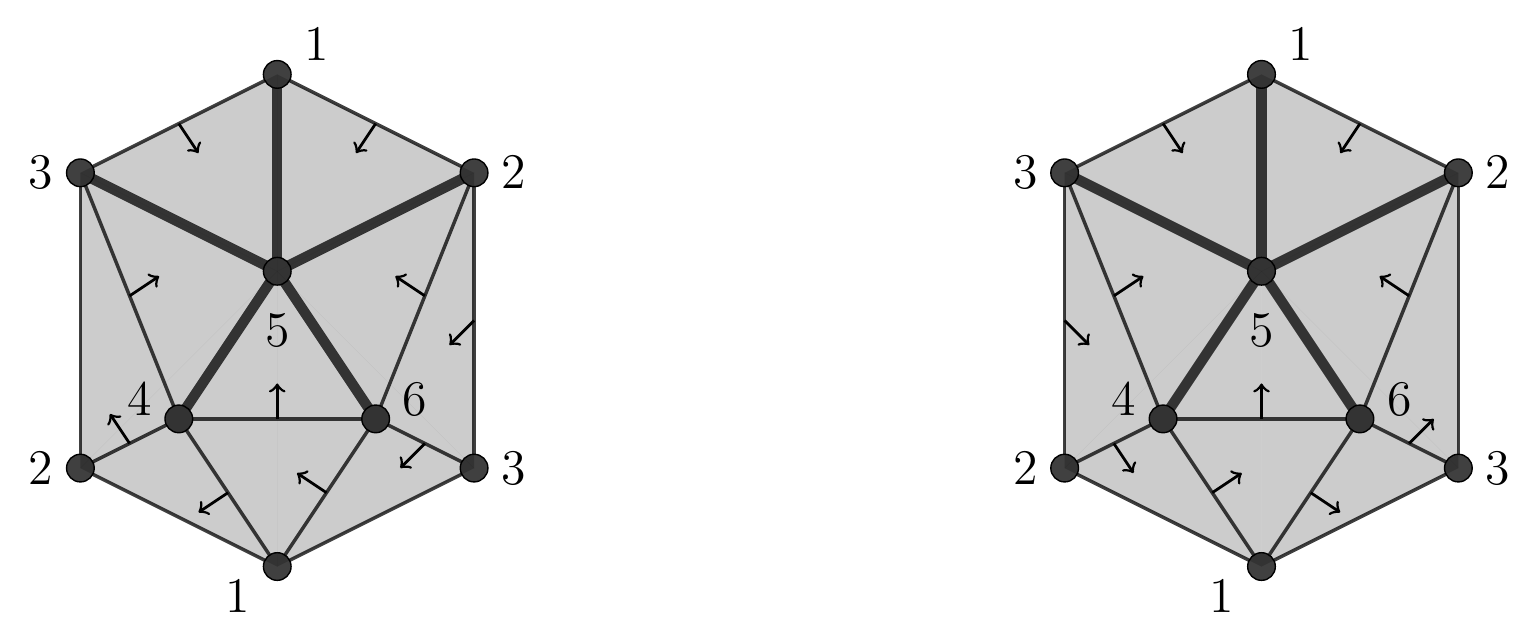}
\end{center}
\caption{Two orientations of the rooted spanning tree of $\Rr P^2$ shown in Fig.~\ref{fig:rooted}~(c). \label{two-orientations}}
\end{figure}

\begin{example} \label{exa:orientations}
The rooted trees (a), (b) and (d) in Figure~\ref{fig:rooted}  each induce a unique orientation.  The rooting of $\Rr P^2$ shown in~(c) induces two distinct orientations: for example, edge $23$ can be paired with either triangle $236$ or triangle $234$.  See Figure~\ref{two-orientations}.  In general, the number of valid orientations is at least $|\HH_{d-1}(Y;\Zz)|$, where $Y$ is the relative complex $(F,R)$, but equality need not hold,
as shown by Bernardi and Klivans~\cite{BK}.
\end{example}

\subsection{Other formulations of higher-dimensional trees} \label{subsection:other-formulations}

The definition of a $d$-forest fits naturally into matroid theory (for a general
reference, see, e.g., \cite{Oxley}).  The {\bf{$d^{th}$ cellular
  matroid}} of a cell complex~$X$ is the matroid $M_d(X)$ represented
over~$\Rr$ by the columns of~$\bd_d(X)$.  This is an extension of the
definition of a \emph{simplicial matroid}; see \cite{CorLin}.
Proposition~\ref{characterize-forest} then asserts that the spanning
$d$-forests of~$X$ are precisely the independent sets in $M_d(X)$, and
the maximal spanning forests are its bases.  For a graph, this
approach would define spanning trees as the bases of the graphical
matroid.

Higher-dimensional trees originate in the work of Bolker~\cite{Bolker}, Kalai~\cite{Kalai} and Adin~\cite{Adin}, who studied spanning trees of skeletons of simplices and joins of 0-dimensional complexes, respectively.
The foregoing definitions, in varying generality and essentially equivalent but for terminological differences, were formulated by Duval, Klivans and Martin \cite{DKM1,DKM2,DKM4}, Lyons~\cite{Lyons}, and Petersson~\cite{Petersson}.  The present definition of a maximal spanning forest of a cell complex $X$ matches the definition of ``spanning tree'' given by Petersson~\cite[Defn.~3.13]{Petersson}.  In addition, Petersson~\cite[Defn.~3.14]{Petersson} defined a pure connected complex to be a tree if its top homology vanishes (i.e., if it is a forest according to Definition~\ref{define-cellular-forest}); by contrast, here Definition~\ref{define-cellular-tree} emphasizes that the higher-dimensional analogue of connectedness is the vanishing of codimension-1 homology.  Lyons~\cite{Lyons} did not use the terms ``tree'' and ``forest,'' instead using matroid terminology.

The literature contains many other generalizations of trees to higher
dimension: see, e.g., \cite{BP,Dewdney,Faridi,HLM,MV,BPT}. In general
these do not have the same topological properties as the cellular trees
considered here, but are defined so as to extend other useful properties of trees.
For instance, Faridi's definition~\cite{Faridi} extends the recursive description
of a tree as a graph formed from a smaller tree by attaching a leaf.  Faridi's trees are simplicial complexes whose dual graphs are trees (such as the first, but not the second, tree in Figure~\ref{fig:trees}); the presence of a leaf is important in the study of their facet ideals.  In another direction, Masbaum and
Vaintrob~\cite{MV} studied 3-uniform hypergraphs (which can be regarded
as pure 2-dimensional simplicial complexes) and found a Pfaffian
matrix-tree theorem enumerating their spanning trees.
\section{Enumerating spanning trees and forests} \label{section:enumeration}

This section will be somewhat more technical and algebraic, with the goal of giving precise statements of the various forms of the cellular matrix-tree theorem that have appeared in the literature and explaining how they are interrelated. 

\subsection{The classical matrix-tree theorem} \label{subsection:classical-MTT}

Let $G=(V,E)$ be a loopless graph with $n=|V|$ and $m=|E|$.  Let $e(v,w)$ denote the number of edges between vertices $i$ and $j$, let $\deg_G(i)$ be the degree of $i$, and let $\bd$ be the signed incidence matrix of $G$ (choosing an arbitrary orientation).
The \defterm{ Laplacian matrix} of $G$ is $L(G)=\bd\bd^T$.
Explicitly, $L(G)=[\ell_{ij}]_{i,j=1}^n$, where
\[\ell_{ij} =
\begin{cases}
\deg_G(i) & \text{ if $i=j$},\\
-e(i,j) & \text{ if $i\neq j$}.
\end{cases}\]

The matrix-tree theorem enumerates spanning trees of a graph using its Laplacian matrix.

\begin{theorem}[Matrix-Tree Theorem] \label{MTT}
Let $G$ be a connected graph, let $\Tree(G)$ be its set of spanning trees, and let $\tau(G)=|\Tree(G)|$.  Then 
\begin{equation} \label{eigenvalue-tree}
\tau(G)=\frac{ \lambda_1 \cdots \lambda_{n-1}}{n}
\end{equation}
where the $\lambda_i$ are the non-zero eigenvalues of $L(G)$.
Alternatively,  for any $i,j$, 
\begin{equation} \label{MTT-reduced-Laplacian}
\tau(G)=\det\tilde L
\end{equation}
where $\tilde L$ is the \defterm{reduced Laplacian} obtained by deleting the $i^{th}$ row and column for some vertex~$i$ (the choice of~$i$ does not matter).
\end{theorem}
The theorem can be proved by expanding the Laplacian or reduced Laplacian using the Binet-Cauchy formula and observing that the non-vanishing terms correspond exactly to spanning trees.  It is also possible to prove~\eqref{MTT-reduced-Laplacian} using the deletion/contraction recurrence for $\tau(G)$.

\begin{corollary}
Let $G$ be any graph (not necessarily connected), let $\Tree(G)$ be its set of maximal spanning forests, and let $\tau(G)=|\Tree(G)|$.  
Suppose that $G$ has nonzero Laplacian eigenvalues $\lambda_1,\dots,\lambda_{n-c}$ and
connected components $G_1,\dots,G_c$ of orders $n_1,\dots,n_c$.  Then
\begin{equation} \label{eigenvalue-forest}
\tau(G) = \frac{\lambda_1 \cdots \lambda_{n-c}}{n_1\cdots n_c} = \det L_S
\end{equation}
where $L_S$ is the $(n-c)\x(n-c)$ matrix obtained by choosing a vertex in each $G_i$ and deleting the corresponding~$c$ rows and columns from $L$.  Moreover, the number of rooted forests of $G$ is simply the product of the nonzero eigenvalues of $L(G)$.
\end{corollary}
Indeed, $L(G)$ is the block-diagonal sum $L(G_1)\oplus\cdots\oplus L(G_c)$, so its spectrum is the multiset union of their spectra.  Meanwhile, a maximal spanning forest of $L(G)$ is just the union of a spanning tree of each component, and the denominator in~\eqref{eigenvalue-forest} is the number of possible choices of roots (since a root consists of a vertex chosen from each component of the graph).

Let $\ww=\{w_e\}$ be a family of commuting indeterminates indexed by the edges of $K_n$.  Let $G$ be a graph on $[n]$ and set $w_e=0$ if $e\not\in E(G)$.
 The \defterm{weighted Laplacian matrix} of $G$ is $L(G;\ww)=[\hat\ell_{i,j}]_{i,j=1}^n$, where
\[\hat\ell_{ij} =
\begin{cases}
\sum_{k\neq i} w_{ik} & \text{ if $i=j$},\\
-w_{ij} & \text{ if $i\neq j$}.
\end{cases}\]

The matrix-tree theorem may be extended to the weighted case, producing polynomial tree enumerators that carry finer combinatorial meaning.  For example, let $x_1,\dots,x_n$ be indeterminates corresponding to the vertices of $G$ and set $w_{ij} = x_ix_j$.  Then  the ``Cayley-Pr\"ufer''  weighted enumeration formula is
\begin{equation} \label{cayley-prufer}
\sum_{T\in\Tree(K_n)} \prod_{i=1}^n x_i^{\deg_T(i)} = x_1\cdots x_n(x_1+\cdots+x_n)^{n-2}.
\end{equation}

\subsection{Extending the matrix-tree theorem to arbitrary cell complexes} \label{subsection:cellular-MTT}

The maximal spanning forests of a cell complex can be enumerated using a higher-dimensional analogue of the matrix-tree theorem.  As described above, this idea originated with Kalai \cite{Kalai} and Adin \cite{Adin} and was extended further by others \cite{CCK,DKM1,DKM2,DKM4,Lyons,Petersson}.

Throughout, let $X$ be a $d$-dimensional pure cell complex. Recall that the {$k^{\rm th}$ boundary map} $\bd_k:\Rr^{X_k}\to\Rr^{X_{k-1}}$ may be regarded as the incidence matrix from $k$-cells to $(k-1)$-cells.  Its adjoint map is the $k^{\rm th}$ \defterm{coboundary map} $\cbd_k$, which may be regarded as the transpose $[\cbd_k]=[\bd_k]^T$.  The $k^{\rm th}$ \defterm{up-down}, \defterm{down-up}, and \defterm{total Laplacians} are respectively 
\[
\Lud_k(X) = \bd_{k+1}\cbd_{k+1},\qquad
\Ldu_k(X) = \cbd_k\bd_k,\qquad
\Ltot_k(X) = {\Lud_k}+\Ldu_k,
\]
all of which may be viewed either as linear operators on $\Rr^{X_k}$, or equivalently as square matrices of size $|X_k|$.  They are all real, diagonalizable, positive semi-definite matrices.  The spectrum of $\Ltot_k(X)$ is the multiset union of the spectra of $\Lud_k(X)$ and $\Ldu_k(X)$ (because the up-down and down-up Laplacians annihilate each other).  Moreover, the spectra
of $\Lud_k(X)$ and $\Ldu_{k+1}(X)$ coincide up to the multiplicity of the zero eigenvalue.  Therefore, the spectra of any of the families $\{\Lud_k(X)\}$, 
$\{\Ldu_k(X)\}$,  $\{\Ltot_k(X)\}$ determines the spectra of the other two.  
The results here are primarily phrased in terms of the up-down Laplacian, so we abbreviate $\Lud_k(X)$ by $L_k(X)$, or simply $L_k$ when the complex $X$ is clear from the context.

Higher-dimensional tree and forest enumeration requires keeping track of torsion.  Accordingly, let $\Tor(A)$ denote the torsion summand of a finitely generated abelian group $A$, and define
\[\tor_k(X)=|\Tor(\HH_k(X;\Zz))|.\]
We extend this notation: if $S$ is any set of cells of $X$, we set $\tor_k(S)=|\Tor(\HH_k(X_S;\Zz))|$, where $X_S$ is the subcomplex generated by $S$.

Let $\Tree_k(X)$ denote the set of maximal spanning $k$-forests in $X$.  As in~\eqref{define-tau} above, the \defterm{$k^{\rm th}$ forest count} of $X$ is
\[
\tau_k(X)=\sum_{T\in\Tree_k(X)}\tor_k(T)^2.
\]
Note that when $\HH_{k-1}(X;\Rr)=0$ (i.e., the spanning forests are spanning trees),
we have $\tor_k(T)=|\Tor(\HH_k(T;\Zz))|=|\HH_k(T;\Zz)|$ for each~$T$.  Moreover, when $X=G$ is a graph, every summand $\tor_0(T)^2$ is 1, so $\tau_1(G)$ is just the number of maximal spanning forests.

\subsection{Weightings} \label{subsection:weightings}

As for graphs, one can study weighted analogues of $\tau(X)$ in order to obtain more refined information about its spanning trees and forests.

Let $\ww=\{w_F\st F\in X\}$ be a family of commuting indeterminates indexed by the faces of~$X$.  The \defterm{weighted forest enumerator} of $X$ is defined as
\begin{equation} \label{define-tau-w}
\tau(X;\ww)=\tau_k(X;\ww) ~\eqdef~ \sum_{T\in\Tree_k(X)} \tor_{k-1}(T)^2\prod_{F\in T} w_F.
\end{equation}
When $X$ is a simplicial complex, we frequently work with variables $v_i$ indexed by the vertices of $X$, and set $w_\sigma=\prod_{i\in \sigma} v_i$ for every face~$\sigma$.  We will refer to this weighting as the \defterm{vertex weighting} for $X$, and denote it by $\vv$ throughout.

Let $D_k$ denote the diagonal matrix with entries $(w_\sigma:\ \sigma\in X_k)$.
As in \cite{CCK,DKM2}, one can define\footnote{This weighted boundary was defined as $\bd_k D_k$ in~\cite{DKM2}.  Here we adopt the equivalent convention of~\cite{CCK} in order to eliminate squares on the indeterminates in weighted enumerators.}
a ``combinatorial'' weighted boundary map $\bd_k D_k^{1/2}$, giving rise to the weighted Laplacian
\begin{equation} \label{combinatorial-weighted-boundary}
L_{k-1}(X;\ww) = \bd_k \,D_k\, \cbd_k.
\end{equation}
Note that the combinatorial weighted boundary maps do not form a chain complex.
Another approach \cite{AD,MMRW,KookLee-WTNMC} is to
define the weighted boundary map  as $D_{k-1}^{-1/2}\bd_k D_k^{1/2}$.
These weighted boundaries  \emph{do} give rise to a chain complex, and to the ``algebraic'' weighted Laplacian:
\begin{equation} \label{algebraic-weighted-boundary}
\Lalg_{k-1}(X;\ww)=D_{k-1}^{-1/2}\,\bd_k\, D_k\, \cbd_k\, D_{k-1}^{-1/2}.
\end{equation}

\subsection{Cellular matrix-forest formulas} \label{subsection:matrix-tree-formulas}

The forest counts of a cell complex~$X$ can be computed from the spectra of its Laplacians.  This relationship appears in several different forms in the literature, all of which share the same underlying principle.  We refer to these results collectively as \emph{cellular matrix-forest theorems}.

The relationship between forest counts and the Laplacian becomes simplest when certain rational homology groups of $X$ vanish.  In \cite{DKM1,DKM2} it was assumed that $X$ was {\bf{$\Rr$-acyclic in positive codimension}} ($\Rr$-APC), i.e., that $\HH_k(X;\Rr)=0$ for all $k<d$. This requirement is not necessary  and so we weaken the hypotheses slightly.  In fact, many complexes of combinatorial interest satisfy a stronger condition: they are \defterm{$\Zz$-acyclic in positive codimension} ($\Zz$-APC), i.e., $\HH_k(X;\Zz)=0$ for all $k<\dim X$. 

For $S\subseteq X_k$, the notation $L_S$ indicates the square
submatrix of $L_k(X)$ restricted to rows and columns indexed by the cells
in $S$.  Define the \defterm{pseudodeterminant} $\pdet(M)$ of any
square matrix $M$ to be the product of its nonzero eigenvalues.

\begin{theorem}[\bf{Cellular-Matrix Forest Theorem}] \label{WCMTT}
Let $X^d$ be a cell complex with $\HH_{d-1}(X;\Rr)=\HH_{d-2}(X;\Rr)=0$.  Let $R\in\Tree_{d-1}(X)$ and $S=X_{d-1}\sm R_{d-1}$. Then
\begin{subequations}
\begin{align}
\tau_d(X;\ww)
&~=~ \frac{\tor_{d-2}(X)^2}{\tor_{d-2}(R)^2}\det\, L_S(X;\ww) \label{WCMTT:eqn2} \\
&~=~ \frac{\padup\tor_{d-2}(X)^2}{\tau_{d-1}(X)} \pdet(L(X;\ww)).\label{WCMTT:eqn1}
\intertext{Furthermore, for a general cell complex $X^d$ with  root $R$ (see Definition~\ref{defn:root}),}
\tau_d(X;\ww)&~=~\frac{\tor_{d-1}(X)^2}{\tor_{d-1}(X,R)^2} \det\, L_S(X;\ww). \label{WCMTT:eqn3}
\end{align}
\end{subequations}
\end{theorem}

These three formulas are respectively~\cite[Thm.~1.4]{DKM1}, \cite[Thm.~2.8]{DKM2}, and~\cite[Prop.~3.5]{DKM4}.\footnote{There is a slight notational error in \cite{DKM4}, where the $L_\Gamma$ in Prop.~3.5 should refer to the reduced Laplacian formed by \emph{removing} the rows and columns corresponding to cells of~$\Gamma$.}
In simple terms,~\eqref{WCMTT:eqn2} and~\eqref{WCMTT:eqn3} state that the forest counts of a complex can be computed from its reduced combinatorial Laplacian, and~\eqref{WCMTT:eqn1} states that the forest counts can be computed via a product of eigenvalues of the combinatorial Laplacian.

When $X$ is $\Zz$-APC, the correction factors in Theorem~\ref{WCMTT} (e.g., $\tor_{d-1}(X)$) are equal to $1$ and therefore may be omitted.

\begin{example}
Let $X$ be the equatorial bipyramid of Example~\ref{example:B}; this complex is $\Zz$-APC.  Let $R$ be the root consisting of all edges that contain vertex $1$.  Then Theorem~\ref{WCMTT} states that the number of spanning trees equals the determinant of the $5\x 5$ reduced Laplacian formed by removing all rows and columns corresponding to edges in $R$; this determinant is $15$.    
\end{example}

The formula~\eqref{WCMTT:eqn1} can be iterated to obtain a closed-form expression for the tree count as an alternating product~\cite[Cor.~2.10]{DKM2}.
For simplicity, we state only the $\Zz$-APC case.
\begin{corollary} \label{cor:alt-product}
Suppose that $\HH_k(X;\Zz)=0$ for all $k<d$.  Then
\begin{equation} \label{alt-product}
\tau_d(X)=\prod_{i=0}^d \pdet(L_{i-1})^{(-1)^{d-i}}.
\end{equation}
\end{corollary}

\begin{example}
Let $X$ be the triangulation of the projective plane as given in Figure~\ref{fig:RP2}.  The pseudodeterminants $\pdet_i=\pdet(L_i)$ of its various Laplacians are:
\[
\pdet_2 = 3^4(3 + \sqrt{5})^3(3- \sqrt{5})^3 = 3^44^3, \qquad
\pdet_1 = 6^5, \qquad
\pdet_0 = 6.
\]
Hence $\tau_2(X)=\pdet_2\pdet_0/\pdet_1=4$.
\end{example}

The right-hand side of Equation~\ref{alt-product} is essentially the R-torsion (or Reidemeister torsion) of~$X$ with respect to the bases of the chain groups given by cells, in the sense of Milnor \cite[\S8]{Milnor}.
Catanzaro, Chernyak and Klein showed more generally~\cite[Thm.~F]{CCK} that the square of the R-torsion is the quantity
\[\prod_{k\geq0}\left(\frac{\covol\HH_k(X;\Zz)\covol B_k(X;\Zz)}{\tor_k(X)^2}\tau_{k+1}(X)\right)^{(-1)^k}\]
where $B_k(X;\Zz)=\im\bd_{k+1}(X)$.

Catanzaro, Chernyak and Klein~\cite{CCK} obtained a closely related determinant formula for $\tau_d(X)$ in terms of the restriction $L_B$ of $L_{d-1}(X)$ to the space $B=B_{d-1}=\im\bd_{d-1}(X)$ (which is invariant under the action of $L_{d-1}$).  For a lattice (discrete subgroup) $\LL\subset\Rr^n$, define its \defterm{covolume}, $\covol\LL$, to be the volume of the parallelotope spanned by any (hence every) set of lattice generators.  Then:

\begin{theorem}[\bf{Cellular Matrix-Forest Theorem, covolume expansion}] \cite[Thms.~C and D]{CCK} \label{MTT-CCK}
\[\tau_d(X;\ww)=\frac{\tor_{d-1}(X)^2}{(\covol B)^2}\det L_B(X;\ww).\]
\end{theorem}

Theorems~\ref{WCMTT} and~\ref{MTT-CCK} are very close in spirit.  Both express the torsion-weighted tree counts of $X$ in terms of determinants of reduced Laplacians.  Theorem~\ref{WCMTT} is useful in combinatorial applications, particularly where there is a natural choice of root $R$ (such as when $X$ is a shifted simplicial complex; see \S\ref{subsection:shifted} below).  Theorem~\ref{MTT-CCK} is less combinatorial but has the advantage of expressing $\tau(X)$ purely in terms of invariants of $X$, with no choice required.

Lyons~\cite{Lyons} gave another general formula for $\tau_k(X)$ (which he denoted by $h_{k-1}(X)$) without the requirements that $\betti_{d-1}(X)=\betti_{d-2}(X)=0$.
For a set of $k$-cells $S\subseteq X_k$, let $\tor'_k(S)$ denote the size of the group $(\ker_\Zz\bd_k)/((\ker_\Zz\bd_k\cap\im_\Qq\bd_{k+1})+\ker_\Zz(\bd_k|_{X_k\sm S}))$,
and define a torsion-weighted enumerator for cobases by
\begin{equation}\label{define-hkX}
h'_k(X)=\sum_{\text{$(k+1)$-cobases }S} \tor_k(X_{k+1}\sm S) \tor'_k(S)^2.
\end{equation}

\begin{theorem}[\bf{Cellular Matrix-Forest Theorem, cobase expansion}] \cite[Prop.~6.2]{Lyons}\label{MTT-Lyons}
Let $S$ be a $(d-1)$-cobase of $X$, i.e., the complement of a $(d-1)$-root $R$.  Then
\begin{subequations}
\begin{align}
\label{Lyons:eqn1}
\tau_d(X)~&=~\frac{\tor_{d-2}(X)^2}{\tor_{d-2}(R)^2\,\tor'_{d-1}(S)^2}\det L_{S}(X)\\
\label{Lyons:eqn2}
~&=~\frac{\padup\tor_{d-2}(X)^2}{h'_{d-2}(X)} \pdet(L(X)).
\end{align}
\end{subequations}
\end{theorem}
When $\HH_{d-1}(X;\Rr)=\HH_{d-2}(X;\Rr)=0$, these formulas reduce to~\eqref{WCMTT:eqn1} and~\eqref{WCMTT:eqn1} respectively.

The algebraic weighting~\eqref{algebraic-weighted-boundary} can be used to obtain a slightly different matrix-forest formula.  Versions of the following result appear in Aalipour, Duval and Martin~\cite{AD}, Kook and Lee~\cite[Thm.~4]{KookLee-WTNMC} and Martin, Maxwell, Reiner and Wilson~\cite[Thm.~5.3]{MMRW}.

\begin{theorem}[\bf{Cellular Matrix-Forest Theorem with algebraic weighting}]\label{WCMTT-AD}
Let $X$ be a cell complex of dimension~$d$ with weighting $\ww$, and suppose that $\HH_{d-1}(X;\Rr)=\HH_{d-2}(X;\Rr)=0$.  Let $\Lalg_{d-1}(X)$ denote the algebraically weighted Laplacian of $X$ as defined in~\eqref{algebraic-weighted-boundary}.  Then
\begin{equation} \label{WCMTT-eqn2}
\tau_d(X;\ww) = \frac{\tor_{d-2}(X)^2}{\tau_{d-1}(X;\ww)}  \pdet(\Lalg_{d-1}(X))
\prod_{\sigma\in X_{d-1}} w_\sigma.
\end{equation}
\end{theorem}
The index $d$ may be replaced with any $k$ for which $\HH_{k-1}(X;\Rr)=\HH_{k-2}(X;\Rr)=0$.
Observe that setting all the weights $w_i$ to 1 recovers~\eqref{WCMTT:eqn1}.

Just as iterating~\eqref{WCMTT:eqn1} gives~\eqref{alt-product}, iterating Theorem~\ref{WCMTT-AD} yields an alternating product formula for weighted tree enumerators.
For simplicity, we state only the $\Zz$-APC case.
\begin{corollary} \label{weighted-alt-prod}
Suppose that $X$ is $\Zz$-APC of dimension~$d$, i.e., $\HH_i(X;\Zz)=0$ for all $i<\dim X$.  Then
\[\tau_d(X;\ww) = \prod_{\substack{\sigma\in X\\ \dim\sigma<d}} (w_\sigma)^{(-1)^{d-\dim\sigma-1}} \prod_{k=-1}^{d-1} \pdet(\Lalg_{k}(X))^{(-1)^{d-k-1}}. \]
\end{corollary}

It is also possible to enumerate rooted forests in terms of products of Laplacian eigenvalues, as observed by Bernardi and Klivans.  

\begin{theorem}[{\bf Rooted Cellular Matrix-Forest Theorem}]~\cite{BK}\label{charpoly}
Let $X$ be a $d$-dimensional cell complex.  Then the characteristic polynomial of the Laplacian $L_d(X)$ gives a generating function for the rooted forests of $X$:
\begin{equation} \label{rooted-CMFT}
\sum_{(F,R)} \tor_{d-1}(X,R)^2 z^{|R|} = \det(L_d(X) + z \cdot \Id).
\end{equation}
\end{theorem}
Here $\Id$ is the identity matrix with rows and columns indexed by~$X_{d-1}$.
The result may also be
extended to include indeterminate weights on the faces of the forest
and root; see \cite[Thm.~18]{BK} and \cite[Prop.~2.1]{MMRW}.
\section{Enumeration of trees for specific complexes} \label{section:specific-enumeration}

We next describe how the theorems of Section~\ref{section:enumeration}
have been used to obtain enumerative results for several families of
cell complexes.  In each case we have provided the full statements of
the tree enumeration formulas but refer the reader to the references throughout for
detailed  treatments of these results.   

\subsection{Skeletons of simplices}\label{subsection:skeletons}

Fix integers $n>d\geq 0$, and let $V=[n]=\{1,2,\dots,n\}$.  The \defterm{complete $d$-dimensional complex on $n$ vertices} is the simplicial complex
\[\Delta=\Delta_{n,d}=\{\sigma\subseteq V\st\dim\sigma\leq d\}.\]
Equivalently, $\Delta$ is the $d$-dimensional skeleton of the simplex on vertices $[n]$.
The complex $\Delta$ is $\Zz$-APC, with top Betti number $\binom{n-1}{d+1}$.

Kalai~\cite[Thm.~1]{Kalai} proved that
\begin{equation} \label{Kalai-formula}
\tau_d(\Delta_{n,d})=n^{\binom{n-2}{d}}
\end{equation}
confirming Bolker's~\cite{Bolker} guess that this formula was correct and generalizing Cayley's formula $\tau(K_n)=n^{n-2}$.  As described in~\S\ref{subsection:intro-enumeration}, Kalai's crucial observation was that the invariant $\tau_d$ is the right way to enumerate higher-dimensional trees.  The actual \emph{number} of spanning trees probably does not have a closed-form expression; it deviates from Kalai's formula first when $n=6$ and $d=2$, since the real projective plane can be triangulated with six vertices; see Example~\ref{K62}.  To obtain formula~\eqref{Kalai-formula}, Kalai
showed that the eigenvalues of the reduced Laplacian are 1 (with multiplicity $\binom{n-2}{d-1}$) and $n$ (with multiplicity $\binom{n-2}{d})$, from which~\eqref{Kalai-formula} follows immediately.

The weighted Cayley-Pr\"ufer formula~\eqref{cayley-prufer} also extends nicely to higher dimensions.  
Kalai~\cite[Thm.~$3'$]{Kalai} proved that
\begin{equation} \label{Kalai-formula-weighted}
\tau_d(\Delta_{n,d};\vv) = (v_1v_2\cdots v_n)^{\binom{n-2}{d-1}} (v_1+v_2+\cdots+v_n)^{\binom{n-2}{d}}
\end{equation}
where $\vv$ is the vertex weighting.  It is an open problem to give a bijective proof of this formula.  It is not clear how to extend the Pr\"ufer code or Joyal's bijection to arbitrary dimension, or how such a bijection should account for trees with torsion.

\subsection{Complete colorful complexes}\label{subsection:complete-colorful}

Let $r$ and $n_1,\dots,n_r$ be positive integers, and let $V_1,\dots,V_r$ be pairwise-disjoint vertex sets with $|V_i|=n_i$.
The \defterm{complete colorful complex} is the $(r-1)$-dimensional complex defined as
\[K_{n_1,\dots,n_r}=\{\sigma\subset V_1\cup\cdots\cup V_r\st \#(\sigma\cap V_i)\leq 1\ \ \forall i\}.\]
The complexes $K_{n_1,\dots,n_r}$ are natural generalizations of complete bipartite graphs, which arise as the case $r=2$.  As described in~\S\ref{subsection:intro-high-dim-trees}, the history of higher-dimensional trees owes its origins to the appearance of complete colorful complexes in the transportation problems considered by Bolker~\cite{Bolker}.  The pure full-dimensional subcomplexes of $K_{n_1,\dots,n_r}$ are precisely the \emph{balanced simplicial complexes}, which are significant in combinatorial commutative algebra; see, e.g.,~\cite[pp.~95--106]{Stanley-CCA}.

As observed by Bolker, the complete colorful complex $K_{n_1,\dots,n_r}$ is $\Zz$-APC, with top Betti number 
$\prod_{i=1}^{r-1} (n_i - 1)$.
Adin \cite[Thm.~1.5]{Adin} gave the following closed formula for their tree counts.
For $D\subseteq[r]$, set
$\sigma_D=\sum_{i\in[r]\sm D} n_i$ and $\pi_D=\prod_{i\in D}(n_i-1)$.
Then for all $0\leq k\leq r-1$,
\begin{equation} \label{Adin-formula}
\tau_k(K_{n_1,\dots,n_r}) = \prod_{\substack{D\subseteq[r]\\ |D|\leq k}} \sigma_D^{\binom{r-2-|D|}{k-|D|}\pi_D}.
\end{equation}
In the special case $k=r-1$, Adin's formula becomes
\begin{equation} \label{Adin-Bolker-formula}
\tau_{r-1}(K_{n_1,\dots,n_r}) = \prod_{i=1}^r n_i^{\left(\prod\limits_{j\neq i}(n_j-1)\right)}.
\end{equation}
Bolker~\cite[Thms.~27,~28]{Bolker} observed that~\eqref{Adin-Bolker-formula} gives the exact count of trees if and only if at most two of the $n_i$'s exceed 2 (because the smallest balanced torsion tree is a 9-vertex triangulation of the real projective plane $\Rr P^2$).

Some additional special cases are worth mentioning.
First, setting $r=2$ in~\eqref{Adin-Bolker-formula} recovers the formula~\eqref{count-Kmn} for the number of spanning trees of a complete bipartite graph.
Second, if $n_i=1$ for all $i$, then $K_{n_1,\dots,n_r}$ is just the simplex on $r$ vertices, so $\pi_D=0$ for all $D$ except $D=\0$, when $\pi_\0=1$, and $\sigma_\0=r$, so~\eqref{Adin-formula} specializes to Kalai's formula~\eqref{Kalai-formula}.
Third, if $n_i=2$ for all $i$, then $K_{n_1,\dots,n_r}$ is the boundary of the $r$-dimensional cross-polytope.  Now $\pi_D=1$ for all $D$ and $\sigma_D=2(r-|D|)$, so~\eqref{Adin-formula} becomes
\[\tau_k(K_{2,2,\dots,2}) = \prod_{d=0}^k (2(r-d))^{\binom{r}{d}\binom{r-2-d}{k-d}}.\]
We will have more to say about this case in \S\ref{subsection:hypercubes} below.

Aalipour, Duval and Martin \cite[Thm.~5.8]{AD} used the weighted matrix-tree theorem (Thm.~\ref{WCMTT-AD}) to prove a vertex-weighted version of Adin's formula~\eqref{Adin-formula}, as follows.  Let $\vv=\{v_{j,t}\st j\in[r],\; t\in[n_j]\}$ be a family of variables indexed by the vertices of $K_{n_1,\dots,n_r}$.
For $j\in[r]$ and $0\leq k\leq r-1$, define
\begin{align*}
p_j &= \prod_{t=1}^{n_j} v_{j,t},&
s_j &= \sum_{t=1}^{n_j} v_{j,t},&
e_{k,j} &= \sum_{\substack{D\subseteq[r]\\j\not\in D,\ |D|\leq k-1}}\!\!\!\!\!\! (-1)^{k-1-|D|}\prod_{t\in D} n_t.
\end{align*}
Then
\begin{equation} \label{Aalipour-Duval-formula}
\tau_k(K_{n_1,\dots,n_r};\vv) = \left(\prod_{j=1}^r p_j^{e_{k,j}}\right)
\prod_{\substack{D\subseteq[r]\\ |D|\leq k}} \left(\sum_{j\in[r]\sm D} s_j\right)^{\binom{r-2-|D|}{k-|D|}\pi_D}.
\end{equation}
Note that setting all $v_{j,t}$ to 1 specializes $p_j$ to 1 and $\sum_{j\not\in D} s_j$ to $\sigma_D$, recovering~\eqref{Adin-formula}.

Complete colorful complexes play a significant role in the theory of
\emph{cyclotomic matroids}, as studied by Martin and
Reiner~\cite{JLM-Cyclotomic}.
In the irreducible case, $M_n$ is isomorphic to the dual of the simplicial
matroid of the complete colorful complex $K(n):=K_{p_1,\dots,p_r}$ where $n$ is the product of distinct primes $p_1, \ldots, p_r$.
That is, the  bases of the matroid are
precisely the complements of the facet sets of spanning trees of
$K(n)$.  Musiker and Reiner~\cite{Musiker-Reiner} extended this work into a beautiful topological
interpretation of the coefficients of the $n^{th}$ cyclotomic
polynomial $\Phi_n$.  In particular, Musiker and Reiner's results
explain why $n=105=3\cdot5\cdot7$ is the first value for which the
coefficients of $\Phi_n$ are not all $\pm1$: it is the smallest $n$
for which the complete colorful complex $K(n)$ contains a spanning
tree with torsion.  In related work, Beck and
Ho\c{s}ten~\cite{Beck-Hosten} studied \emph{cyclotomic polytopes},
which are dual to the generalized transportation polytopes considered
by Bolker, and observed that patterns for their growth series break
down for similar reasons involving torsion trees.

\subsection{Shifted simplicial complexes}\label{subsection:shifted}
Let $\sigma=\{i_0<\dots<i_s\}$ and $\rho=\{j_0<\cdots<j_r\}$ be sets of positive integers, listed in increasing order.  Write $\sigma\preceq\rho$ if $\sigma\subseteq\rho$, or if $s=r$ and $i_k\leq j_k$ for every $k$.  This relation is a partial order, called the \defterm{componentwise order}.  
A simplicial complex~$\Delta^d$ is called \defterm{shifted} if its vertices can be labeled $1,\dots,n$ so that in each dimension its faces form an order ideal with respect to this order: that is, if $\rho\in\Delta$ and $\sigma\preceq\rho$, then $\sigma\in\Delta$ as well.  Shifted complexes are an important class of extremal complexes; every simplicial complex can be transformed into a shifted complex in such a way that many combinatorial properties are preserved while the overall structure is simplified.  A signal success of the study of shifted complexes was the classification of all $f$-vector/Betti-number pairs for all simplicial complexes, by Bj\"{o}rner and Kalai \cite{Bjorner-Kalai}.
Shifted complexes of dimension~1 are known as \emph{threshold graphs} and have generated much study in their own right; see \cite{MahadevPeled}.

The homology groups of a shifted complex~$\Delta$  are free abelian; the Betti number $\betti_k(\Delta)$ is the number of facets of $\Delta$ of dimension~$k$ not containing vertex~1.
To see this, observe that contracting the subcomplex $S=\Star_\Delta(1)$ produces a wedge of spheres, one for each facet not containing vertex~1.
In particular, if $\Delta$ is pure then it is $\Zz$-APC, and $S$ is a $\Zz$-acyclic spanning tree.  For the remainder of this section, we assume that $\Delta$ is pure of dimension~$d$ on vertex set $[n]$.

Duval and Reiner~\cite{Duval-Reiner} showed that shifted simplicial complexes are Laplacian integral.  Specifically, let $\deg_\Delta(i)$ denote the number of facets containing vertex~$i$, so that $(\deg_\Delta(1),\ldots,\deg_\Delta(n))$ is a partition.  Then the conjugate (transpose) partition is precisely the spectrum of 
the top-dimensional Laplacian $L_{d-1}(\Delta)$, 
up to the multiplicity of the zero eigenvalue.  This generalizes the corresponding result for threshold graphs, first observed by Hammer and Kelmans \cite[Thm.~5.3]{HK} and stated in this form by Merris~\cite{Merris}.

Weighted enumerators for spanning trees of threshold graphs were calculated by Martin and Reiner \cite[\S6]{JLM-Factorization},
and can be deduced as a special case of a more general bijective result of Remmel and Williamson \cite[Thm.~2.4]{RW}.
Weighted tree enumeration for pure shifted complexes was studied by Duval, Klivans and Martin \cite{DKM1}, as we now summarize.
Introduce doubly-indexed indeterminates $\ww=\{w_{i,j} \st 0\leq i\leq d,\ 1\leq j\leq n\}$ and for each face $\sigma=\{a_0<\cdots<a_k\}$, set
\(w_\sigma = \prod_{i=0}^k w_{i,a_i}\).
In other words, $w_\sigma$ records not just which vertices occur in $\sigma$, but their order within $\sigma$.  Theorem 1.5 of~\cite{DKM1} states that the eigenvalues
of the weighted Laplacian $L_k(\Delta;\ww)$ are polynomials in the weights, so that setting all weights to~1 gives an explicit proof of Laplacian integrality for shifted complexes.  Specifically, the eigenvalues are in bijection with \defterm{critical pairs} of $\Delta$: pairs $(\sigma,\rho)$ of $(k+1)$-sets of vertices such that $\sigma\in\Delta$, $\rho\not\in\Delta$, and $\sigma$ is covered by $\rho$ in the componentwise order.  

Let $\Gamma=\{\rho\in\Delta \st 1\not\in\rho\}$ and $\Lambda=\{\rho\in \Gamma \st \rho\cup\{1\}\in\Delta\}$ denote respectively the deletion and link of vertex~1 (which themselves are shifted complexes).
Observe that for every critical pair $(\sigma,\rho)=((a_0,\dots,a_d),(b_0,\dots,b_d))$, there is a unique index~$j$ for which $b_j=a_j+1$
(and $b_i=a_i$ for all~$i\neq j$).  Define the \defterm{signature} of~$(\sigma,\rho)$ to be $(S,T)$, where $S=\{a_0,\dots,a_{j-1}\}$ and $T=\{1,2,\dots, a_j\}$.  Then \cite[Thm.~1.6]{DKM1}
\begin{equation}\label{shifted-formula-fine}
\tau_d(\Delta;\ww) 
= \left( \prod_{F\in \Lambda_{d-1}} w_{F\cup\{1\}} \right)
\left( \prod_{(S,T)\in\Sig(\Gamma)} \frac{\sum\limits_{j \in T\cup\{1\}} w_{S \cup j}}{w_{S\cup\{1\}}} \right)
\end{equation}
where $\Sig(\Gamma)$ denotes the set of signatures of critical pairs of~$\Gamma$.  
The case $d=1$ (i.e., $\Delta$ is a threshold graph) is equivalent to~\cite[Thm.~4]{JLM-Factorization}.  Specializing $w_{i,j}\mapsto v_j$ for every $i,j$ to obtain the vertex weighting $\vv$ (called the ``coarse weighting'' in \cite{DKM1}) produces the formula
\begin{equation}\label{shifted-formula-coarse}
\tau_d(\Delta;\vv)
~=~ v_1^{|\Lambda_{d-1}|} \left(\prod_{j=2}^n v_j^{\deg_\Lambda(j)}\right)
\left( \prod_{(S,T)\in\Sig(\Gamma)} v_1^{-1}\sum_{j \in T\cup\{1\}} v_j \right).
\end{equation}

\subsection{Color-shifted complexes}\label{subsection:color-shifted}

As in \S\ref{subsection:complete-colorful}, let $K=K_{n_1,\dots,n_r}$ be the complete colorful complex with vertex set $V_1\cup\cdots\cup V_r$, where $V_i=\{v_{i,1},\dots,v_{i,n_i}\}$.  A subcomplex $\Delta\subseteq K$ is called \defterm{color-shifted} if it satisfies the following property:
\begin{equation} \label{color-shifted-defn}
\text{if } \sigma=\{v_{1,a_1},v_{2,a_2},\dots,v_{r,a_r}\}\in\Delta \text{ and } 1\leq b<a_i, \text{ then } \sigma\sm\{v_{i,a_i}\}\cup\{v_{i,b}\}\in\Delta.
\end{equation}
Color-shifted complexes were introduced by Babson and Novik~\cite{Babson-Novik}.  Structurally, color-shifted complexes are to complete colorful complexes as shifted complexes are to simplices.  On the other hand, color-shifted complexes (even in dimension~1) are not in general Laplacian integral.

Let $G$ be a color-shifted complex with $r=2$, i.e., a graph.  Assume that $G$ is connected, and let $n=n_1$, $m=n_2$.
For $1\leq i\leq n$, let $\lambda_i=\max\{j: v_{1,i}v_{2,j}\in E(G)\}$.  Then $\lambda=(\lambda_1,\dots,\lambda_n)$ is a partition, and it determines $G$ completely: the edges of $G$ correspond to the boxes in the Ferrers diagram of $\lambda$, with
$V_1$ and $V_2$ corresponding to its rows and columns respectively.  For this reason, 1-dimensional color-shifted complexes were studied under the name of \emph{Ferrers graphs} by Ehrenborg and van~Willigenburg, who found an elegant formula for their weighted tree enumerators~\cite[Thm.~2.1]{EvW}.
Let $G=G_\lambda$ be the Ferrers graph corresponding to a partition~$\lambda$,
and weight the vertices of $V_1$ and $V_2$ by $\xx=\{x_1,\dots,x_n\}$ and $\yy=\{y_1,\dots,y_m\}$ respectively.
Then
\begin{equation} \label{Ferrers-graphs}
\tau(G_\lambda;\xx,\yy) = x_1\cdots x_n y_1\cdots y_m \prod_{p=1}^n(y_1+\cdots+y_{\lambda_p-1})
\prod_{q=1}^m(x_1+\cdots+x_{\tilde\lambda_q-1})
\end{equation}
where $\tilde\lambda$ is the conjugate partition to~$\lambda$.  The proof of this formula in~\cite{EvW} is recursive, proceeding by
calculating $\tau(G_\lambda;\xx,\yy)/\tau(G_{\lambda'};\xx,\yy)$, where $\lambda'$ is a partition obtained from~$\lambda$ by removing a corner of the Ferrers diagram.  Another proof, expressing $\tau(G_\lambda;\xx,\yy)$ in terms of the weighted tree enumerator of a related threshold graph, was given in \cite[\S10.3]{DKM1}.

For a color-shifted complex~$\Delta$ with $r$ colors, let $\vv=\{v_{q,i}\}$ be the vertex weighting, where $q\in[r]$ indexes the colors and $i$ indexes the vertices of each color.  
Let $\Delta'$ be the subcomplex of $\Delta$ whose faces are $\{\rho \in \Delta\st v_{q,1} \not\in \rho\ \text{for any $q$}\}$.
Note that every $(r-2)$-face $\sigma\in\Delta$ is missing a unique color $m=m(\sigma)$; set
$k(\sigma) = \max\{j\st \sigma\cup\{x_{m,i}\}\in\Delta\}$.   Aalipour and Duval conjecture~\cite{AD-colorshifted} that 
\begin{equation} \label{Aalipour-Duval}
\tau_{r-1}(\Delta;\vv) = \prod_{q,i} x_{q,i}^{e(q,i)}
\prod_{\sigma\in\Delta'_{r-2}}
(x_{m,1}+\cdots+x_{m,k(\sigma)})
\end{equation}
where $e(q,i)=\#\{\sigma\in\Delta_{r-1} \st x_{q,i}\in\sigma \text{ and } x_{q',1}\in\sigma \text{ for some } q'\neq q\}$.
Aalipour and Duval have proven this formula in the case $r=3$ (where a color-shifted complex corresponds to a plane partition).  It reduces to \eqref{Ferrers-graphs} when $r=2$, and to the weighted tree enumerator formula~\eqref{Aalipour-Duval-formula} when $\Delta$ is a complete colorful complex.

Equation~\eqref{Aalipour-Duval} resembles formulas~\eqref{shifted-formula-fine} and~\eqref{shifted-formula-coarse} for weighted tree enumerators of shifted complexes.  Specifically, the first factor on the right-hand side is a product of monomials indexed by faces \emph{with} vertices numbered 1, and the second factor is a product of sums indexed by faces \emph{without} vertices numbered 1.

\subsection{Matroid complexes}\label{subsection:matroid}

Let $M$ be a matroid on ground set $V$.  The family $\II(M)$ of independent sets of $M$ forms a simplicial complex called a \defterm{matroid complex}.  These complexes can be characterized in many ways; perhaps the simplest is that a simplicial complex is a matroid complex if and only if every induced subcomplex is pure.  Matroid complexes are  Cohen-Macaulay and $\Zz$-APC.  (See Oxley~\cite{Oxley} for a general reference on matroids, and Stanley~\cite[Chap.~III]{Stanley-CCA} for facts about matroid complexes.)

Every matroid complex is Laplacian integral.  This fact was proved first by Kook, Reiner and Stanton \cite{KRS}, and subsequently generalized by Denham~\cite{Denham}, who found an explicit eigenbasis for the Laplacian.
An expression for $\tau(\II(M))$ can be extracted from Denham's result (with difficulty) by applying~\eqref{WCMTT:eqn1} or~\eqref{alt-product}.
One application of the Kook-Reiner-Stanton formula, given in \cite[\S6.1]{DKM2}, is an alternate proof of Adin's formula \eqref{Adin-formula}, because the complete colorful complex $K_{n_1,\dots,n_r}$ is the independence complex of a direct sum of rank-one uniform matroids.

Kook and Lee~\cite{KookLee-WTNMC,KookLee-FSTNMC} gave a simpler expression for weighted tree enumerators of matroid complexes.
Let $\alpha(M) =T_M(0,1)$, where $T_M$ is the Tutte polynomial of~$M$ (see~\cite{BryOx})
and $\beta(M) = (-1)^{r(M)}\sum_{A\subseteq E}(-1)^{|A|}r(A)$ (the Crapo beta invariant~\cite{Crapo-beta}, where $r$ is the rank function; see also~\cite[pp.~447--448]{Stanley-HA}).
Then
\[\tau(\II(M);\vv) = \prod_{e\in E} v_e^{|\Tree(M/e)|-\alpha(M/e)} \prod_{\text{flats } F} \left(\sum_{e\in E\sm F} v_e\right)^{\alpha(F)\beta(M/F)}\]
\cite[Thm.~9]{KookLee-WTNMC}, where $\vv$ is the vertex weighting.  Consequently,
\begin{equation} \label{KookLee}
\tau(\II(M)) = \prod_{\text{flats } F} (|E\sm F|)^{\alpha(F)\beta(M/F)}.
\end{equation}

\subsection{Hypercubes} \label{subsection:hypercubes}
The \defterm{hypercube} $Q_n$ is the space $[0,1]^n\subseteq\Rr^n$,
naturally regarded as a cell complex with cells $\{\{0\},\{1\},[0,1]\}^n$.  Note that the hypercube is not a simplicial complex, but it has a tight combinatorial structure that enables enumeration of its spanning trees.

The hypercube is contractible, hence $\Zz$-APC.  The eigenvalues of $L_0(Q_n,r)$ are $0,2,\dots,2j,\dots,2n$, with multiplicities $\binom{n}{j}$; this can be shown by directly constructing an eigenbasis as in \cite[pp.~61--62]{EC2}.
Applying formula~\eqref{eigenvalue-tree} of the classical matrix-tree theorem yields the formula
\begin{equation} \label{count-Qn}
\tau_1(Q_n)=2^{-n}\prod_{A\subseteq[n]}2|A|=2^{2^n-n-1}\prod_{k=2}^n k^{\binom{n}{k}}.
\end{equation}
This formula was discovered (in greater generality) by Cvetkovi\'c~\cite[(4.10)]{Cvetkovic}, using spectral techniques; see also \cite[pp.~75, 218]{SpectraBook}.
The first purely combinatorial proof was recently found by Bernardi~\cite{Bernardi}.
A generalization to higher dimension was obtained by Duval, Klivans and Martin in \cite[Thm.~3.4]{DKM2}, who showed that the nonzero eigenvalues of $L_k(Q_n)$ are $2(k+1),\dots,2(j+1),\dots,2n$, with multiplicities $\binom{j-1}{k} \binom{n}{j}$.  Applying the alternating-product formula~\eqref{alt-product} gives the tree count
\begin{equation} \label{tau-k-Qn}
\tau_k(Q_n) ~=~ \prod_{j=k+1}^n (2j)^{\binom{n}{j}\binom{j-2}{k-1}}.
\end{equation}

Weighted tree enumeration of hypercubes proceeds similarly.
Introduce indeterminates $\{q_i,x_i,y_i\st i\in[n]\}$ and weight each face 
$f=(f_1,\dots,f_n) \in \{\{0\},\{1\},[0,1]\}^n$ 
by the monomial
\[w_f=\left(\prod_{i:\;f_i=[0,1]}q_i\right)\left(\prod_{i:\;f_i=0}x_i\right)\left(\prod_{i:\;f_i=1}y_i\right).\]
Now consider the algebraically weighted Laplacian $\Lalg_k(Q_n;\ww)$ (see~\eqref{algebraic-weighted-boundary}) for $k\geq0$.
By~\cite[Thm.~4.2]{DKM2}, its nonzero eigenvalues are the rational functions $u(A)=\sum_{i\in A} q_i/x_i+q_i/y_i$, indexed by subsets $A\subseteq[n]$ with $|A|>k$, with multiplicities $\binom{|A|-1}{k}$.
Applying the weighted alternating product formula\footnote{One needs the convention that $w_\0=1$, so that $\Lalg_{-1}(Q_n)$ is just the $1\x1$ matrix with entry $\sum_{v\in V(Q_n)}w_v=\prod_{j=1}^n(x_j+y_j)$.}
(Corollary~\ref{weighted-alt-prod}) yields
\begin{align}
\tau_k(Q_n;\ww	)
~&=~ (q_1\cdots q_n)^{\sum_{i=k-1}^{n-1}\binom{n-1}{i}\binom{i-1}{k-2}}
\prod_{\substack{A\subseteq[n]\\ |A|>k}}
\left(u(A)\prod_{j\in A} x_jy_j\right)^{\binom{|A|-2}{k-1}}
\label{weighted-tau-k-Qn}
\end{align}
as conjectured in \cite[Conj.~4.3]{DKM2}.  For a complete proof of this formula, see~\cite{AD}.  The case $k=1$ gives the weighted enumerator for spanning trees of cube graphs~\cite[Thm.~3]{JLM-Factorization}.

\subsection{Duality and self-duality} \label{subsection:duality}

Let $X,Y$ be cell complexes of dimension~$d$.  We say that $X$ and $Y$ are \defterm{dual} if there is an inclusion-reversing bijection $\sigma\mapsto\sigma^*$ between their face posets that induces a duality between their cellular chain complexes.
On the level of chain complexes, this means that $\bd_{X,k}=\cbd_{Y,d-k+1}$ for every $k$, so\footnote{In~\cite{DKM2}, it was erroneously stated that $\Lud_k(X)=\Lud_{d-k}(Y)$.  In fact, duality interchanges up-down and down-up Laplacians.} $\Lud_k(X)=\Ldu_{d-k}(Y)$.
Indeed, $T\subseteq X$ is a spanning $k$-forest in $X$ if and only if $T^\vee=\{\sigma^*\st\sigma\in X\sm T\}$ is a spanning $k$-forest in $Y$ \cite[Prop.~6.1]{DKM2}.
Duality carries over to the weighted setting: given a weighting~$\ww$ on $X$, define the dual weighting $\ww^*$ on $Y$ by $w^*(\sigma^*)=w(\sigma)^{-1}$.  Then $\Lud_k(X;\ww)=\Ldu_k(Y;\ww^*)$, and likewise for the algebraically weighted Laplacians (see \S\ref{subsection:weightings}).  By the weighted cellular matrix-tree theorem (Thm.~\ref{WCMTT}), $\tau_k(X;\ww)=\tau_{d-k}(Y;\ww^*)$ for all $k$.
As an application, the $(n-1)$-skeleton of the hypercube $Q_n$ is dual to the complete colorful complex $K_{2,2,\dots,2}$ (with $n$ 2's), so the formulas~\eqref{tau-k-Qn} and~\eqref{weighted-tau-k-Qn}
for spanning tree enumerators of $Q_n$ are dual to special cases of~\eqref{Adin-formula} and~\eqref{Aalipour-Duval-formula} respectively.

A \defterm{central reflex} is an embedding of a self-dual planar graph $G$ on the sphere $S^2$ so that the isomorphism $G\to G^*$ (where $G$ denotes the planar dual) is realized by the antipodal map on $\Ss^2$.
Tutte~\cite{Tutte-self-dual} proved that if $G$ is a central reflex, then $\tau(G)$ is the square of the number of self-dual spanning trees of $G$.   Maxwell~\cite{Maxwell} generalized Tutte's result to higher dimension.  Say that a matroid on ground set $E=\{e_1,\dots,e_n,\tilde e_1,\dots,\tilde e_n\}$
is \defterm{involutively self-dual} if it can be represented by an $n\x 2n$ integer matrix $M$ of the form $[I|N]$, where $I$ is the $n\x n$ identity matrix and $N$ is skew-symmetric (so that the fixed-point free involution $E\to E$ given by $e_i\longleftrightarrow\tilde e_i$ gives an isomorphism $M\to M^*$).  If $\Delta$ is the independence complex of such a matroid (so $\dim\Delta = d = n-1$), then \cite[Thm.~1.3]{Maxwell}
\begin{equation} \label{Maxwell-eqn}
\tau_d(\Delta) = \sum_{T\in\Tree_d(\Delta)}\tor_{d-1}(T)^2 = \left( \sum_{\text{self-dual }T\in\Tree_d(\Delta)}\tor_{d-1}(T)\right)^2.
\end{equation}

As an application, let $k$ be a positive integer and let $\Delta$ be a simplex on vertex set $[2k+2]$.  Then $\tau_k(\Delta_{2k+2,k}) = (2k+2)^{\binom{2k}{k}}$ by Kalai's formula~\eqref{Kalai-formula}.  The number $\binom{2k}{k}$ is even for all $k>0$, so this tree enumerator is a perfect square.  The simplicial matroid of the $k$-skeleton of $\Delta$ is involutively self-dual via the map $\sigma\mapsto[2k+2]\sm\sigma$
 so this perfect square phenomenon is explained by~\eqref{Maxwell-eqn}.

On the level of weighted tree enumerators, the formula for $\tau_k(\Delta_{2k+2,k};\ww)$ given by~\eqref{Kalai-formula-weighted} is a perfect square when $n=2k+2$ and $k$ is odd (since in that case $\binom{n}{k-1}$ must be even).
Maxwell~\cite[Thm.~1.8]{Maxwell} proved that this perfect square phenomenon is combinatorial, in the sense that
the square root of $\tau_k(\Delta_{2k+2,k};\ww^2)$ has the combinatorial interpretation one would expect from~\eqref{Maxwell-eqn}:
\[\sum_{T\in\Tree_k(\Delta)}\tor_{k-1}(T)^2 w_T^2
~=~ \left( \sum_{\text{self-dual }T\in\Tree_k(\Delta)}\tor_{k-1}(T) w_T\right)^2.\]
The proof identifies the left-hand side with the determinant of a skew-symmetric matrix, and the weighted enumerator for self-dual trees with its Pfaffian.

Self-dual balls of \emph{even} dimension were considered by Martin, Maxwell, Reiner and Wilson in \cite{MMRW}.
If $X$ is a self-dual cellular ball of dimension $d=2k$ and $\bd=\bd_k(X)$, then $\pdet(L_{k-1}(X))=\tau_{k-1}(X)\tau_k(X)=\tau_k(X)^2$.  More generally, consider a weighting~$\ww$ such that $w_\sigma=w_{\sigma^*}$, where $\sigma^*$ is the dual cell  to~$\sigma$.  Then $\pdet(\Lalg_{k-1}(X;\ww)) = \tau_k(X;\ww)^2$, a consequence of \cite[Cor.~6.10]{MMRW}.  If in addition the self-duality on~$X$ is induced by an antipodal map, then it can be shown that~$X$ has an orientation giving rise to a symmetric or skew-symmetric boundary matrix~$\bd$ (depending on the parity of~$k$).  In either case $\pdet(L_{k-1})=\pdet(\bd\cbd)=\pdet(\bd)^2$, so in fact $\tau_k(X)=\tau_{k-1}(X)=\pdet(\bd)$ and more generally $\tau_k(X;\ww)=\pdet(D_k\bd)$ by~\cite[Thm.~7.4]{MMRW}.
\section{Connections, further directions and open problems} \label{section:open-problems}

\subsection{Critical groups of cell complexes} \label{subsection:critical}

The theory of \emph{critical groups} emerges from a discrete dynamical system that has been studied in mathematics, physics, probability, and even economics in several closely related forms, including the \emph{abelian sandpile model}, the \emph{chip-firing game}, and the \emph{dollar game}.  Original sources include \cite{Bak,Biggs,BLS,Dhar,Lorenzini2}.  See \cite{LevinePropp} for a short accessible introduction to the subject, and \cite{Merino-chip-firing} for an overview.

The critical group can be defined purely algebraically from the Laplacian of a graph.  In higher dimension, the algebraic definitions have natural analogues using the tools developed earlier in the chapter.

Let $G=(V,E)$ be a graph and let $L=L(G)$ be the graph Laplacian.
The \defterm{critical group} $K(G)$ is the quotient group
\[K(G) = \Tor( \Zz^V/\im\, L ) = \ker \bd_0 / \im\, L \, \cong \, \Zz^{|V|-1}/\im\,\tilde L\]
where $\tilde L$ is the reduced Laplacian (see Theorem~\ref{MTT}).
The isomorphism is due to Biggs~\cite[Thm.~4.2]{Biggs}.
Regarding a vector $\cc\in\Zz^V$ as a configuration in which $c_i$ is the number of grains of sand (or chips, or dollars) on vertex~$v_i$, the elements of the critical group can be thought of as equivalence classes of configurations modulo a dispersive Laplacian action.  A special system of coset representatives, the \emph{critical configurations}, encodes the possible long-term stable behaviors of the system.

By the matrix-tree theorem, the order of the critical group, and hence the number of distinct long-term configurations, both equal the number of spanning trees of $G$.  There are deep connections between critical configurations and spanning trees, as well as the combinatorics of parking functions, $G$-parking functions, and hyperplane arrangements; see, e.g.,~\cite{Biggs-Winkler,CoriLB,ChPy,BCT,Primer,BS, GK}.

Accordingly, the \defterm{$i^{th}$ critical group} $K_i(X)$ of a cell complex $X$ is defined as
$$
K_i(X) = \Tor(\Zz^{X_i}/\im\, L_i) = \ker \bd_{i}/\im\, L_i \, \cong \,  \Zz^{|S|} / \im\, L_S
$$
where $L_S$ is the reduced Laplacian obtained by removing a torsion-free $i$-tree, provided one exists.  The isomorphism was proven by Duval, Klivans and Martin~\cite[Thm.~3.4]{DKM3} in the simplicial case; the proof carries over to the more general cellular setting.  In analogy to the graph case, the size of $K_i(X)$ is the forest count $\tau_{i+1}(X)$. 

These higher critical groups model discrete dynamical systems that generalize the sandpile model~\cite[\S5]{DKM3}.  A state of the system described by $K_i(X)$ consists of a number for each $i$-cell; this quantity can be regarded as ``charge'' ($i=0$), ``current'' ($i=1$), ``circulation'' ($i=2$), etc.  The Laplacian action redistributes this quantity between $i$-cells that lie in a common $(i+1)$-cell.  When $i=0$, this system is the sandpile model described briefly above.  When $i=1$, each state assigns a current to each edge, and redistribution takes the form of circulation around 2-cells.
Because the configurations are contained in $\ker \bd_i$, the
$i$-dimensional flow is \defterm{conservative}: the sum of oriented
flow at every $i$-dimensional face containing a fixed
$(i-1)$-dimensional face is zero.  When $i=1$, the
oriented flow on the edges neither accumulates nor depletes
at any vertex.  When $i=2$, the face circulation at each
edge adds up to zero.  

The model can  be regarded as having a sink consisting
of an $i$-dimensional spanning tree. By
specifying flow values on all $i$-faces that are {\em not} in the
sink, the values on the sink needed to make the flow conservative are
uniquely defined; this is why the configurations can be described by the cokernel of the reduced Laplacian.
Guzm\'an and Klivans~\cite{GK2} have recently identified meaningful systems of representatives of the higher critical group which reflect the long term behavior of the Laplacian action.

\subsection{Cut and Flow Spaces} \label{subsection:cutflow}

The critical group of a cell complex is closely related to its cuts and flows, which generalize the corresponding objects in graph theory.  In a connected graph $G$, an \defterm{(edge) cut} is a collection
of edges whose removal disconnects the graph, and a
\defterm{bond} is a minimal cut.  In matroid terms, a bond is a cocircuit: a minimal set of edges that meets every basis.  The support of every row of the incidence matrix $\bd$ is a cut, namely the set of all edges incident to a particular vertex, and more generally the \defterm{cut space} is defined as $\im\cbd$, that is, the row space of $\bd$.

Meanwhile, a \defterm{flow} on $G$ is an assignment of a ``current'' to each edge so that the total flow into each vertex equals the total flow out (called a conservative flow above); that is, a flow is an element of the \defterm{flow space} $\ker\bd=(\im\cbd)^\perp$.  In matroid terms, the minimal supports of flows are circuits of the graphic matroid (i.e., simple cycles in $G$).  The \defterm{cut lattice} $\C(G)$ and \defterm{flow lattice} $\F(G)$ are defined in the same way as the cut and flow spaces, except regarding $\bd$ as a map of $\Zz$-modules rather than vector spaces.

Accordingly, for a cell complex $X^d$, the \defterm{cut space}, \defterm{cut lattice}, \defterm{flow space}, and \defterm{flow lattice} of $X$ in dimension~$k$ are defined as
\begin{align*}
\Cut_k(X) &= \im_{\Rr} \cbd_k(X), & \C_k(X) &= \im_\Zz \cbd_k(X),\\
\Flow_k(X) &= \ker_{\Rr} \bd_k(X), & \F_k(X) &= \ker_\Zz \bd_k(X).
\end{align*}
An $k$-dimensional \defterm{bond} of $X$ is a minimal set of $k$-faces that support a non-zero vector in $\Cut_k(X)$, and a $k$-dimensional \defterm{circuit} of $X$ is a minimal set of $k$-faces that support a non-zero vector of $\Flow_k(X)$.  (These are exactly the cocircuits and circuits, respectively, of the cellular matroid of~$X$; see \S\ref{subsection:trees-and-matroids} and \cite{SuWagner}.)  Flows and cuts can also be defined homologically: e.g., vectors in $\Flow_d(X)$ are just elements of $\HH_d(X;\Rr)$, while vectors in $\Cut_d(X)$ are sets of facets whose deletion increases codimension-1 homology.

Given a spanning tree $T$ of a connected graph~$G$, the (signed) characteristic vectors of the fundamental bonds and fundamental cycles of $T$ form bases of the cut and flow spaces respectively, and in fact they are $\Zz$-module bases for the cut and flow lattices.  This is a standard fact of algebraic graph theory~\cite[Chap.~14]{Godsil-Royle} that carries over well to general cell complexes~\cite[\S\S4--6]{DKM4}.  Predictably, it is necessary to account for torsion in higher dimension.  The entries of the characteristic vectors of the fundamental bonds and circuits of a cellular spanning tree $T\subseteq X$
can be interpreted as cardinalities of certain torsion groups \cite[Thms.~4.11 and~5.3]{DKM4}, and in order to form $\Zz$-module bases of $\C_d(X)$ or $\F_d(X)$, the tree $T$ must be $\Zz$-acyclic.

For an integral lattice $\LL$ (that is, a subgroup of $\Zz^n$), the \defterm{dual lattice} is $\LL^\sharp=\{\ww\in\Rr^n\st\langle\ww,\vv\rangle\in\Zz\ \ \forall\vv\in\LL\}$, where $\langle\cdot,\cdot\rangle$ is the standard inner product.  Thus $\LL\subseteq\LL^\sharp$, and $\LL^\sharp/\LL$ is called the \defterm{discriminant group} (or determinantal group) of~$\LL$.  Bacher, de~la~Harpe and Nagnibeda~\cite{BHN} proved that the groups
$K(G)$, $\C^\sharp/\C$, $\F^\sharp/\F$, and $\Zz^{E(G)}/(\C\oplus\F)$
are mutually isomorphic.  (See also~\cite[Chap.~14]{Godsil-Royle}.)
The relationship between the cut and flow lattices, their discriminant groups, and the critical group of $X$ is more complicated in arbitrary dimension.  The isomorphisms are replaced with short exact sequences
\begin{align*}
&0 ~\to~ \Zz^n / (\C\oplus\F) ~\to~  \hspace{.5em} K(X) \isom \C^\sharp/\C ~\to~ E ~\to~ 0,\\
0 ~\to~ &E ~\to~ \Zz^n/(\C\oplus\F)~\to~ K^*(X) \isom \F^\sharp/\F ~\to~ 0 & 
\end{align*}
where $E$ is the ``error term'' $\Tor(\HH_{d-1}(X;\Zz))$~\cite[Thms.~7.6 and~7.7]{DKM4},
and $K^*(X)$ is the \defterm{cocritical group}, defined roughly by dualizing the definition of the critical group~\cite[Defn.~7.4]{DKM4}. 
In particular, the groups $K(X)$, $K^*(X)$, $\C^\sharp/\C$, $\F^\sharp/\F$, and $\Zz^{E(G)}/(\C\oplus\F)$ are all isomorphic if and only if $\HH_{d-1}(X;\Zz)$ is torsion-free.  

Bacher, de la Harpe and Nagnibeda's work in \cite{BHN} arose from viewing a graph as a combinatorial analogue of a Riemann surface, where the critical group appears as the discrete analogue of the Picard group of divisors, or the Jacobian group of holomorphic forms.  This point of view also informs the graphical Riemann-Roch theorem of Baker and Norine~\cite{Baker-Norine} and subsequent work connected to tropical geometry~\cite{Luo,HMY,HKN}.
  An open problem posed in~\cite[\S6]{DKM3} is to understand how the combinatorics of critical and cocritical groups on cell complexes parallels the geometric theory of group invariants such as Chow groups on algebraic varieties.

\subsection{Trees and simplicial complex decomposition} \label{subsection:decomp}

Recall from Definitions~\ref{defn:root} and~\ref{defn:orientation} that a \defterm{directed rooted spanning forest} of a cell complex $X^d$ is a rooted spanning forest $(F,R)$, together with an orientation $\OO$ that pairs each $(d-1)$-face $\rho$ outside the root~$R$ with a facet $\phi\in F$ that contains~$\rho$.  As seen in Example~\ref{exa:orientations}, rooted forests of dimension $\geq2$ can admit more than one orientation: as shown by Bernardi and Klivans~\cite{BK}, the number of valid orientations is at least $|\HH_{d-1}(F,R)|$.  Equality need not hold: for instance, rooted triangulations of the (torsion-free) dunce cap can admit multiple orientations.  In an appropriate linear-algebraic sense, these excess orientations cancel in pairs to the size of the torsion subgroup.  It is an open problem to give a combinatorial explanation of this cancellation; see~\cite{BK}.

Iterating the construction of a directed rooted spanning forest gives
a discrete (although not gradient) vector field in the sense of
discrete Morse theory~\cite{Morse1,Morse2}.  To be precise, let
$(F,R,\OO)$ be a directed rooted $d$-forest of $X^d$.  Then $R$ is a $(d-1)$-dimensional subcomplex of $X^d$.  Choose a directed rooted forest of $R$
which further pairs faces of $X^d$, this time of dimensions $(d-2)$ and $(d-1)$.  
Continuing this process, the vector field is then given by the union of the pairings, and the unmatched cells are enumerated by the Betti numbers of~$X$.  For example, consider the directed spanning 2-tree of $\Rr P^2$ shown in (a) or (b) of Figure~\ref{two-orientations}.  This is the 1-dimensional root, and the arrows indicate the bijection $\OO:F\to X_1\sm R$.  In the next iteration, choose a single vertex as the 0-dimensional root, and orient all of the edges of $R$ away from the root to match them with nonroot vertices.  Finally, match the empty face with the root vertex.

These recursive orientations arise in theorems and conjectures about decompositions of simplicial complexes, with applications to their $f$- and $h$-vectors.  In particular, certain results about decompositions of simplicial complexes can be usefully restated in the language of simplicial forests and orientations, as we now describe.

Suppose that $\Delta$ is a simplicial complex that is acyclic over at least one coefficient ring (it does not matter which), so that the $f$-polynomial $f(\Delta,x)=\sum_{\sigma\in\Delta} x^{|\sigma|}$ is divisible by $x-1$.  Stanley explained this algebraic fact combinatorially by showing that the face poset of~$\Delta$ can be decomposed as $\Delta'\dju\Omega$, where $\Delta'$ is a subcomplex and there is a bijection $\eta:\Delta'\to\Omega$ such that every $\sigma\in\Delta'$ is a maximal proper subface of $\eta(\sigma)$ \cite[Thm.~1.2]{Stanley-Decomp}.
Duval~\cite[Thm.~1.1]{Duval-Decomp} showed more generally that every simplicial complex $\Delta$ admits a decomposition $B\dju\Delta'\dju\Omega$, where $B$ contains $\betti_k(\Delta)$ faces for each $k$ and $\Delta',\Omega$ are as just described.  Duval's result (which in fact holds for all regular cell complexes) led to a new proof of Bj\"orner and Kalai's characterization of $f$-vector/Betti-number pairs of arbitrary simplicial complexes~\cite{Bjorner-Kalai}, without using algebraic shifting.
In Stanley's and Duval's constructions, the set $\Omega$ is a \defterm{flag of forests}: an order filter in the face poset of $\Delta$ such that for every $k$, the set $\Omega_k=\{\sigma\in\Omega\st\dim\sigma=k\}$ is the set of facets of a spanning $k$-forest of~$\Delta$.  

More generally, Stanley conjectured~\cite[Conj.~2.4]{Stanley-Decomp} that if $\Delta$ is \defterm{$k$-uply acyclic}, i.e., if the link of every face of dimension $\leq k-2$ is acyclic, then its face poset can be decomposed into Boolean algebras of rank~$k$.  The $k=1$ case is Stanley's result cited above.  The general ce remains open, although Duval and Zhang proved a weakening of Stanley's conjecture in~\cite{Duval-Zhang}.  Stanley introduced this conjecture in the same context as  the \emph{partitionability conjecture}, which stated that every Cohen-Macaulay simplicial complex $\Delta$ can be partitioned into Boolean intervals whose tops are its facets, with ranks given by the $h$-vector.
In light of Remark~\ref{ex:z-acyclic}, there are implications
\begin{equation*}
\parbox[c]{3.5cm}{every link in $X$ has a\\$\Zz$-acyclic spanning tree} \implies \, \,  \parbox{2.0cm}{$X$ is Cohen-\\Macaulay} \implies \, \, \, \, \parbox{2.7cm}{every link in $X$ has\\a spanning tree}
\end{equation*}
because by Reisner's criterion~\cite{Reisner}, a simplicial complex is Cohen-Macaulay if the link of every face is $\Zz$-APC.

The present authors envisioned simplicial spanning trees as an inroad to the partitionability conjecture, culminating in its disproof in~\cite{DGKM}.  Meanwhile, in a different direction, Katth\"an~\cite{Katthan} had discovered that the \emph{depth conjecture} (proposed by Stanley~\cite{Stanley-LDE}; see also~\cite{WhatIs,Herzog-survey}) could be reduced to the study of certain lattices associated to simplicial spanning trees.  While the disproof of the partitionability conjecture implies that the depth conjecture is also false in general, Katth\"an's methods reveal a deeper connection between trees and combinatorial commutative algebra: e.g., he showed that spanning trees of simplex skeletons (which he called \emph{stoss complexes}) arise as the Scarf complexes of certain extremal lattices~\cite[Thm.~4.3]{Katthan}.

\subsection{Trees and Matroids} \label{subsection:trees-and-matroids}

As mentioned in \S\ref{subsection:other-formulations}, the collection of all spanning
forests of a cell complex $X$ forms the collection of bases of a matroid,
the \emph{cellular matroid} $M(X)$, which generalizes the simplicial matroids
of~\cite{CorLin} and coincides with the well-known \emph{graphic matroid} in dimension~$1$.  Matroid structures such as cocircuits are useful in studying cellular trees (as in \S\ref{subsection:cutflow}), so it is natural to ask what other matroid structures and invariants mean in the context of trees.

Krushkal and Renardy~\cite{KR} studied a polynomial invariant on cell complexes that generalizes the Tutte polynomial of the cellular matroid.  (For a general reference on the Tutte polynomial, see~\cite{BryOx}.)
Let $X$ be a cell complex whose geometric realization is an orientable manifold of even dimension. Its Krushkal-Renardy invariant $P_X(x,y,a,b)$ is a generating function for spanning $d$-dimensional subcomplexes $Y\subseteq X$, in which the exponents of $x,y$ record combinatorial information about the cellular matroid and those of $a,b$ record topological information about the embedding of~$Y$ in~$X$.  Setting $a=b=1$ recovers the Tutte polynomial of $M(X)$, and $P$ satisfies a duality relation $P_X(x,y,a,b)=P_{X^*}(y,x,b,a)$, where $X^*$ is the 
dual cell complex \cite[Thm.~14]{KR}, generalizing the usual Tutte duality relation $T_M(x,y)=T_{M^*}(y,x)$ on matroids.

Evaluating the Tutte polynomial $T_M(x,y)$ of a matroid $M$ at $(x,y)=(1,1)$ gives the number of bases of~$M$.  If $M$ is the graphic matroid of a graph $G$ then $T_M(1,1)=\tau(G)$, but this equality does not hold for most cellular matroids because of torsion.  On the other hand, torsion can be modeled combinatorially via the theory of \emph{arithmetic matroids}, introduced by d'Adderio and Moci \cite{Moci1,Moci2,Moci3} and further developed
by Fink and Moci~\cite{FinkMoci}.  An arithmetic matroid is a matroid equipped with
a \emph{multiplicity function} $m$ that models the number of torsion elements in an abelian group.
The arithmetic Tutte polynomial studied in \cite{Moci1,Moci2} records this multiplicity. For a cellular matroid $M(X^d)$, the multiplicity of a set $A$ of facets is $|\Tor(\HH_{d-1}(A\cup X_{\leq d-1};\Zz))|^2$, and
evaluating the arithmetic Tutte polynomial at
$(x,y)=(1,1)$ gives the torsion-weighted tree count of $X$.
Bajo, Burdick and Chmutov~\cite{BBC} similarly modified the Krushkal-Renardy polynomial to include a torsion factor, again obtaining the
torsion-weighted tree count as a suitable evaluation of their polynomial invariant.

The matroid structure for trees suggests that graph algorithms could carry over well to the general setting of cell complexes.  Of particular interest is the problem of sampling random trees, proposed by Russell Lyons \cite{Lyons} and Igor Pak.  In the graph case, there exist efficient algorithms for tree sampling. The best known is Wilson'€™s loop-erased random walk algorithm~\cite{Wilson}, which generates a uniformly distributed random spanning tree. Extending Wilson'€™s algorithm to general cell complexes is likely to be hard, as there is no known efficient algorithm to sample random bases of an arbitrary matroid.   (It is possible for some classes, including equitable matroids; see~\cite{Mayhew}; Gorodezky and Pak~\cite{Igors} generalized Wilson's algorithm to a restricted class of hypergraphs that do not coincide with cellular forests.)  One could look for an algorithm that samples cellular trees either uniformly, or according to a distribution that includes torsion information.
Indeed, this problem informs Lyons' motivation in~\cite{Lyons}, in which he defined determinantal probability measures on cellular trees based on the Laplacian.

\begin{acknowledgement}
We thank our colleagues Amy Wagler and Jennifer Wagner for valuable suggestions regarding exposition and clarity.
\end{acknowledgement}

\bibliographystyle{amsalpha}
\bibliography{biblio}

\def\soft#1{\leavevmode\setbox0=\hbox{h}\dimen7=\ht0\advance \dimen7
  by-1ex\relax\if t#1\relax\rlap{\raise.6\dimen7
  \hbox{\kern.3ex\char'47}}#1\relax\else\if T#1\relax
  \rlap{\raise.5\dimen7\hbox{\kern1.3ex\char'47}}#1\relax \else\if
  d#1\relax\rlap{\raise.5\dimen7\hbox{\kern.9ex \char'47}}#1\relax\else\if
  D#1\relax\rlap{\raise.5\dimen7 \hbox{\kern1.4ex\char'47}}#1\relax\else\if
  l#1\relax \rlap{\raise.5\dimen7\hbox{\kern.4ex\char'47}}#1\relax \else\if
  L#1\relax\rlap{\raise.5\dimen7\hbox{\kern.7ex
  \char'47}}#1\relax\else\message{accent \string\soft \space #1 not
  defined!}#1\relax\fi\fi\fi\fi\fi\fi}
\providecommand{\bysame}{\leavevmode\hbox to3em{\hrulefill}\thinspace}
\providecommand{\MR}{\relax\ifhmode\unskip\space\fi MR }
\providecommand{\MRhref}[2]{%
  \href{http://www.ams.org/mathscinet-getitem?mr=#1}{#2}
}
\providecommand{\href}[2]{#2}
\begin{thebibliography}{MMRW15}

\bibitem[AD]{AD-colorshifted}
Ghodratollah Aalipour and Art~M. Duval, \emph{Weighted spanning tree
  enumerators of color-shifted complexes}, in preparation.

\bibitem[Adi92]{Adin}
Ron~M. Adin, \emph{Counting colorful multi-dimensional trees}, Combinatorica
  \textbf{12} (1992), no.~3, 247--260. \MR{1195888 (93j:05036)}

\bibitem[ADM]{AD}
Ghodratollah Aalipour, Art~M. Duval, and Jeremy~L. Martin, \emph{Weighted
  spanning tree enumerators of complete colorful complexes}, in preparation.

\bibitem[Aus60]{Austin}
T.~L. Austin, \emph{The enumeration of point labelled chromatic graphs and
  trees}, Canad. J. Math. \textbf{12} (1960), 535--545. \MR{0139544 (25
  \#2976)}

\bibitem[BBC14]{BBC}
Carlos Bajo, Bradley Burdick, and Sergei Chmutov, \emph{On the
  {T}utte-{K}rushkal-{R}enardy polynomial for cell complexes}, J. Combin.
  Theory Ser. A \textbf{123} (2014), 186--201. \MR{3157807}

\bibitem[BBGM14]{BBGM}
Matthias Beck, Felix Breuer, Logan Godkin, and Jeremy~L. Martin,
  \emph{Enumerating colorings, tensions and flows in cell complexes}, J.
  Combin. Theory Ser. A \textbf{122} (2014), 82--106. \MR{3127679}

\bibitem[BCT10]{BCT}
Brian Benson, Deeparnab Chakrabarty, and Prasad Tetali, \emph{{$G$}-parking
  functions, acyclic orientations and spanning trees}, Discrete Math.
  \textbf{310} (2010), no.~8, 1340--1353. \MR{2592488 (2011i:05152)}

\bibitem[BdlHN97]{BHN}
Roland Bacher, Pierre de~la Harpe, and Tatiana Nagnibeda, \emph{The lattice of
  integral flows and the lattice of integral cuts on a finite graph}, Bull.
  Soc. Math. France \textbf{125} (1997), no.~2, 167--198. \MR{1478029
  (99c:05111)}

\bibitem[Ber12]{Bernardi}
Olivier Bernardi, \emph{On the spanning trees of the hypercube and other
  products of graphs}, Electron. J. Combin. \textbf{19} (2012), no.~4, Paper
  51, 16. \MR{3007186}

\bibitem[BH06]{Beck-Hosten}
Matthias Beck and Serkan Ho{\c{s}}ten, \emph{Cyclotomic polytopes and growth
  series of cyclotomic lattices}, Math. Res. Lett. \textbf{13} (2006), no.~4,
  607--622. \MR{2250495 (2007h:52018)}

\bibitem[Big99]{Biggs}
N.~L. Biggs, \emph{Chip-firing and the critical group of a graph}, J. Algebraic
  Combin. \textbf{9} (1999), no.~1, 25--45. \MR{1676732 (2000h:05103)}

\bibitem[BK]{BK}
Olivier Bernardi and Caroline~J. Klivans, \emph{Directed rooted forests in
  higher dimension}, in preparation.

\bibitem[BK88]{Bjorner-Kalai}
Anders Bj{\"o}rner and Gil Kalai, \emph{An extended {E}uler-{P}oincar{\'e}
  theorem}, Acta Math. \textbf{161} (1988), no.~3-4, 279--303. \MR{971798
  (89m:52009)}

\bibitem[BLS91]{BLS}
Anders Bj{\"o}rner, L{\'a}szl{\'o} Lov{\'a}sz, and Peter~W. Shor,
  \emph{Chip-firing games on graphs}, European J. Combin. \textbf{12} (1991),
  no.~4, 283--291. \MR{1120415 (92g:90193)}

\bibitem[BN06]{Babson-Novik}
Eric Babson and Isabella Novik, \emph{Face numbers and nongeneric initial
  ideals}, Electron. J. Combin. \textbf{11} (2004/06), no.~2, Research Paper
  25, 23 pp. (electronic). \MR{2195431 (2007c:05202)}

\bibitem[BN07]{Baker-Norine}
Matthew Baker and Serguei Norine, \emph{Riemann-{R}och and {A}bel-{J}acobi
  theory on a finite graph}, Adv. Math. \textbf{215} (2007), no.~2, 766--788.
  \MR{2355607 (2008m:05167)}

\bibitem[BO92]{BryOx}
Thomas Brylawski and James Oxley, \emph{The {T}utte polynomial and its
  applications}, Matroid applications, Encyclopedia Math. Appl., vol.~40,
  Cambridge Univ. Press, Cambridge, 1992, pp.~123--225. \MR{1165543
  (93k:05060)}

\bibitem[Bol76]{Bolker}
Ethan~D. Bolker, \emph{Simplicial geometry and transportation polytopes},
  Trans. Amer. Math. Soc. \textbf{217} (1976), 121--142. \MR{0411983 (54
  \#112)}

\bibitem[Bor60]{Borchardt}
C.~W. Borchardt, \emph{{\"U}ber eine {I}nterpolationsformel f{\"u}r eine {A}rt
  {S}ymmetrischer {F}unctionen und {\"u}ber {D}eren {A}nwendung}, Math. Abh.
  der Akademie der Wissenschaften zu Berlin (1860), 1--20.

\bibitem[BP71]{BP}
L.~W. Beineke and R.~E. Pippert, \emph{Properties and characterizations of
  {$k$}-trees}, Mathematika \textbf{18} (1971), 141--151. \MR{0288046 (44
  \#5244)}

\bibitem[BPT15]{BPT}
Yurii Burman, Andrey Ploskonosov, and Anastasia Trofimova, \emph{Matrix-tree
  theorems and discrete path integration}, Linear Algebra Appl. \textbf{466}
  (2015), 64--82. \MR{3278240}

\bibitem[BS13]{BS}
Matthew Baker and Farbod Shokrieh, \emph{Chip-firing games, potential theory on
  graphs, and spanning trees}, J. Combin. Theory Ser. A \textbf{120} (2013),
  no.~1, 164--182. \MR{2971705}

\bibitem[BTW88]{Bak}
Per Bak, Chao Tang, and Kurt Wiesenfeld, \emph{Self-organized criticality},
  Phys. Rev. A (3) \textbf{38} (1988), no.~1, 364--374. \MR{949160 (89g:58126)}

\bibitem[BW97]{Biggs-Winkler}
Norman Biggs and Peter Winkler, \emph{Chip-firing and the chromatic
  polynomial}, Tech. Report LSE-CDAM-97-03, London School of Economics, Center
  for Discrete and Applicable Mathematics, 1997.

\bibitem[CCK12]{CCK}
Michael~J. Catanzaro, Vladimir~Y. Chernyak, and John~R. Klein,
  \emph{Kirchhoff's theorems in higher dimensions and {Reidemeister} torsion},
  Preprint, \href{http://arxiv.org/abs/1206.6783}{arXiv:1206.6783}, 2012.

\bibitem[CDS80]{SpectraBook}
Drago{\v{s}}~M. Cvetkovi{\'c}, Michael Doob, and Horst Sachs, \emph{Spectra of
  graphs}, Pure and Applied Mathematics, vol.~87, Academic Press, Inc.
  [Harcourt Brace Jovanovich, Publishers], New York-London, 1980, Theory and
  application. \MR{572262 (81i:05054)}

\bibitem[Chu97]{Chung-book}
Fan R.~K. Chung, \emph{Spectral graph theory}, CBMS Regional Conference Series
  in Mathematics, vol.~92, Published for the Conference Board of the
  Mathematical Sciences, Washington, DC; by the American Mathematical Society,
  Providence, RI, 1997. \MR{1421568 (97k:58183)}

\bibitem[CL87]{CorLin}
Raul Cordovil and Bernt Lindstr{\"o}m, \emph{Simplicial matroids},
  Combinatorial geometries, Encyclopedia Math. Appl., vol.~29, Cambridge Univ.
  Press, Cambridge, 1987, pp.~98--113. \MR{921070}

\bibitem[CLB03]{CoriLB}
Robert Cori and Yvan Le~Borgne, \emph{The sand-pile model and {T}utte
  polynomials}, Adv. in Appl. Math. \textbf{30} (2003), no.~1-2, 44--52, Formal
  power series and algebraic combinatorics (Scottsdale, AZ, 2001). \MR{1979782
  (2004d:05095)}

\bibitem[CP05]{ChPy}
Denis Chebikin and Pavlo Pylyavskyy, \emph{A family of bijections between
  {$G$}-parking functions and spanning trees}, J. Combin. Theory Ser. A
  \textbf{110} (2005), no.~1, 31--41. \MR{2128964 (2005m:05010)}

\bibitem[Cra67]{Crapo-beta}
Henry~H. Crapo, \emph{A higher invariant for matroids}, J. Combinatorial Theory
  \textbf{2} (1967), 406--417. \MR{0215744 (35 \#6579)}

\bibitem[Cve71]{Cvetkovic}
Drago{\v{s}}~M. Cvetkovi{\'c}, \emph{The spectral method for determining the
  number of trees}, Publ. Inst. Math. (Beograd) (N.S.) \textbf{11(25)} (1971),
  135--141. \MR{0309772 (46 \#8877)}

\bibitem[Den01]{Denham}
Graham Denham, \emph{The combinatorial {L}aplacian of the {T}utte complex}, J.
  Algebra \textbf{242} (2001), no.~1, 160--175. \MR{1844702 (2002h:05039)}

\bibitem[Dew74]{Dewdney}
A.~K. Dewdney, \emph{Higher-dimensional tree structures}, J. Combinatorial
  Theory Ser. B \textbf{17} (1974), 160--169. \MR{0369115 (51 \#5351)}

\bibitem[DGKM15]{DGKM}
Art~M. Duval, Bennet Goeckner, Caroline~J. Klivans, and Jeremy~L. Martin,
  \emph{A non-partitionable cohen-macaulay simplicial complex}, preprint,
  \href{http://www.arXiv.org/1504.04279}{arXiv:1504.04279}, 2015.

\bibitem[Dha90]{Dhar}
Deepak Dhar, \emph{Self-organized critical state of sandpile automaton models},
  Phys. Rev. Lett. \textbf{64} (1990), no.~14, 1613--1616. \MR{1044086
  (90m:82053)}

\bibitem[DKM09]{DKM1}
Art~M. Duval, Caroline~J. Klivans, and Jeremy~L. Martin, \emph{Simplicial
  matrix-tree theorems}, Trans. Amer. Math. Soc. \textbf{361} (2009), no.~11,
  6073--6114. \MR{2529925 (2011a:05385)}

\bibitem[DKM11]{DKM2}
\bysame, \emph{Cellular spanning trees and {L}aplacians of cubical complexes},
  Adv. in Appl. Math. \textbf{46} (2011), no.~1-4, 247--274. \MR{2794024
  (2012e:05182)}

\bibitem[DKM13]{DKM3}
\bysame, \emph{Critical groups of simplicial complexes}, Ann. Comb. \textbf{17}
  (2013), no.~1, 53--70. \MR{3027573}

\bibitem[DKM15]{DKM4}
\bysame, \emph{Cuts and flows of cell complexes}, J. Algebraic Combin.
  \textbf{41} (2015), 969--999.

\bibitem[DM12]{Moci1}
Michele D'Adderio and Luca Moci, \emph{Ehrhart polynomial and arithmetic
  {T}utte polynomial}, European J. Combin. \textbf{33} (2012), no.~7,
  1479--1483. \MR{2923464}

\bibitem[DM13a]{Moci3}
\bysame, \emph{Arithmetic matroids, the {T}utte polynomial and toric
  arrangements}, Adv. Math. \textbf{232} (2013), 335--367. \MR{2989987}

\bibitem[DM13b]{Moci2}
\bysame, \emph{Graph colorings, flows and arithmetic {T}utte polynomial}, J.
  Combin. Theory Ser. A \textbf{120} (2013), no.~1, 11--27. \MR{2971693}

\bibitem[DP76]{Dodziuk}
J.~Dodziuk and V.~K. Patodi, \emph{Riemannian structures and triangulations of
  manifolds}, J. Indian Math. Soc. (N.S.) \textbf{40} (1976), no.~1-4, 1--52
  (1977). \MR{0488179 (58 \#7742)}

\bibitem[DR02]{Duval-Reiner}
Art~M. Duval and Victor Reiner, \emph{Shifted simplicial complexes are
  {L}aplacian integral}, Trans. Amer. Math. Soc. \textbf{354} (2002), no.~11,
  4313--4344 (electronic). \MR{1926878 (2003j:15017)}

\bibitem[Duv94]{Duval-Decomp}
Art~M. Duval, \emph{A combinatorial decomposition of simplicial complexes},
  Israel J. Math. \textbf{87} (1994), no.~1-3, 77--87. \MR{1286816 (96e:52023)}

\bibitem[DZ01]{Duval-Zhang}
Art~M. Duval and Ping Zhang, \emph{Iterated homology and decompositions of
  simplicial complexes}, Israel J. Math. \textbf{121} (2001), 313--331.
  \MR{1818393 (2003a:52013)}

\bibitem[Eck45]{Eckmann}
Beno Eckmann, \emph{Harmonische {F}unktionen und {R}andwertaufgaben in einem
  {K}omplex}, Comment. Math. Helv. \textbf{17} (1945), 240--255. \MR{0013318
  (7,138f)}

\bibitem[EvW04]{EvW}
Richard Ehrenborg and Stephanie van Willigenburg, \emph{Enumerative properties
  of {F}errers graphs}, Discrete Comput. Geom. \textbf{32} (2004), no.~4,
  481--492. \MR{2096744 (2005j:05076)}

\bibitem[Far02]{Faridi}
Sara Faridi, \emph{The facet ideal of a simplicial complex}, Manuscripta Math.
  \textbf{109} (2002), no.~2, 159--174. \MR{1935027 (2003k:13027)}

\bibitem[FH98]{Friedman-Hanlon}
Joel Friedman and Phil Hanlon, \emph{On the {B}etti numbers of chessboard
  complexes}, J. Algebraic Combin. \textbf{8} (1998), no.~2, 193--203.
  \MR{1648484 (2000c:05155)}

\bibitem[FM12]{FinkMoci}
Alex Fink and Luca Moci, \emph{Matroids over a ring}, J.\ Eur.\ Math.\ Soc., to
  appear; \href{http://www.arXiv.org/1209.6571}{arXiv:1209.6571}, 2012.

\bibitem[For98]{Morse1}
Robin Forman, \emph{Morse theory for cell complexes}, Adv. Math. \textbf{134}
  (1998), no.~1, 90--145. \MR{1612391 (99b:57050)}

\bibitem[For02]{Morse2}
\bysame, \emph{A user's guide to discrete {M}orse theory}, S\'em. Lothar.
  Combin. \textbf{48} (2002), Art.\ B48c, 35. \MR{1939695 (2003j:57040)}

\bibitem[Fri98]{Friedman}
J.~Friedman, \emph{Computing {B}etti numbers via combinatorial {L}aplacians},
  Algorithmica \textbf{21} (1998), no.~4, 331--346. \MR{1622290 (99c:52022)}

\bibitem[FS58]{Fiedler-Sedlacek}
Miroslav Fiedler and Ji{\v{r}}{\'{\i}} Sedl{\'a}{\v{c}}ek, \emph{{\"U}ber
  {W}urzelbasen von gerichteten {G}raphen}, {\v C}asopis P{\v e}st. Mat.
  \textbf{83} (1958), 214--225. \MR{0097071 (20 \#3551)}

\bibitem[GK15a]{GK}
Johnny Guzm{\'a}n and Caroline Klivans, \emph{Chip-firing and energy
  minimization on {M}-matrices}, J. Combin. Theory Ser. A \textbf{132} (2015),
  14--31. \MR{3311336}

\bibitem[GK15b]{GK2}
\bysame, \emph{Chip-firing on general invertible matrices}, in preparation,
  2015.

\bibitem[GP14]{Igors}
Igor Gorodezky and Igor Pak, \emph{Generalized loop-erased random walks and
  approximate reachability}, Random Structures Algorithms \textbf{44} (2014),
  no.~2, 201--223. \MR{3158629}

\bibitem[GR01]{Godsil-Royle}
Chris Godsil and Gordon Royle, \emph{Algebraic graph theory}, Graduate Texts in
  Mathematics, vol. 207, Springer-Verlag, New York, 2001. \MR{1829620
  (2002f:05002)}

\bibitem[Hat02]{Hatcher}
Allen Hatcher, \emph{Algebraic topology}, Cambridge University Press,
  Cambridge, 2002. \MR{1867354 (2002k:55001)}

\bibitem[Her13]{Herzog-survey}
J{\"u}rgen Herzog, \emph{A survey on {S}tanley depth}, Monomial ideals,
  computations and applications, Lecture Notes in Math., vol. 2083, Springer,
  Heidelberg, 2013, pp.~3--45. \MR{3184118}

\bibitem[HK96]{HK}
P.~L. Hammer and A.~K. Kelmans, \emph{Laplacian spectra and spanning trees of
  threshold graphs}, Discrete Appl. Math. \textbf{65} (1996), no.~1-3,
  255--273, First International Colloquium on Graphs and Optimization (GOI),
  1992 (Grimentz). \MR{1380078 (97d:05205)}

\bibitem[HKN13]{HKN}
Jan Hladk{\'y}, Daniel Kr{\'a}{\soft{l}}, and Serguei Norine, \emph{Rank of
  divisors on tropical curves}, J. Combin. Theory Ser. A \textbf{120} (2013),
  no.~7, 1521--1538. \MR{3092681}

\bibitem[HLM99]{HLM}
J{\"u}rgen Herzog and Enzo~Maria Li~Marzi, \emph{Bounds for the {B}etti numbers
  of shellable simplicial complexes and polytopes}, Commutative algebra and
  algebraic geometry ({F}errara), Lecture Notes in Pure and Appl. Math., vol.
  206, Dekker, New York, 1999, pp.~157--167. \MR{1702104 (2001b:13017)}

\bibitem[HMY12]{HMY}
Christian Haase, Gregg Musiker, and Josephine Yu, \emph{Linear systems on
  tropical curves}, Math. Z. \textbf{270} (2012), no.~3-4, 1111--1140.
  \MR{2892941}

\bibitem[Joy81]{Joyal}
Andr{\'e} Joyal, \emph{Une th{\'e}orie combinatoire des s{\'e}ries formelles},
  Adv. in Math. \textbf{42} (1981), no.~1, 1--82. \MR{633783}

\bibitem[Kal83]{Kalai}
Gil Kalai, \emph{Enumeration of {${\bf Q}$}-acyclic simplicial complexes},
  Israel J. Math. \textbf{45} (1983), no.~4, 337--351. \MR{720308 (85a:55006)}

\bibitem[Kat14]{Katthan}
Lukas Katth\"an, \emph{Stanley depth and simplicial spanning trees}, J.\
  Algebraic Combin., to appear;
  \href{http://www.arXiv.org/1410.3666}{arXiv:1410.3666}, 2014.

\bibitem[Kir47]{Kirchhoff}
G.~Kirchhoff, \emph{{\"U}œber die {A}ufl{\"o}sung der {G}leichungen, auf welche
  man bei der {U}ntersuchung der linearen {V}erteilung galvanischer
  {S}tr{\"o}me gef{\"u}hrt wird}, Ann. Phys. Chem. \textbf{72} (1847),
  497--508.

\bibitem[KL15a]{KookLee-FSTNMC}
Woong Kook and Kang-Ju Lee, \emph{A formula for simplicial tree numbers of
  matroid complexes}, 2015, preprint.

\bibitem[KL15b]{KookLee-WTNMC}
\bysame, \emph{Weighted tree-numbers of matroid complexes}, 27th
  {I}nternational {C}onference on {F}ormal {P}ower {S}eries and {A}lgebraic
  {C}ombinatorics ({FPSAC} 2015), Discrete Math. Theor. Comput. Sci. Proc., AS,
  Assoc. Discrete Math. Theor. Comput. Sci., Nancy, 2015, to appear.

\bibitem[KR14]{KR}
Vyacheslav Krushkal and David Renardy, \emph{A polynomial invariant and duality
  for triangulations}, Electron. J. Combin. \textbf{21} (2014), no.~3, Paper
  3.42, 22. \MR{3262279}

\bibitem[KRS00]{KRS}
W.~Kook, V.~Reiner, and D.~Stanton, \emph{Combinatorial {L}aplacians of matroid
  complexes}, J. Amer. Math. Soc. \textbf{13} (2000), no.~1, 129--148.
  \MR{1697094 (2001e:05028)}

\bibitem[KW68]{Klee-Witzgall}
Victor Klee and Christoph Witzgall, \emph{Facets and vertices of transportation
  polytopes}, Mathematics of the {Decision} {Sciences}, {Part} {I} ({Seminar},
  {Stanford}, {Calif}., 1967), Amer. Math. Soc., Providence, R.I., 1968,
  pp.~257--282. \MR{0235832}

\bibitem[Lor91]{Lorenzini2}
Dino~J. Lorenzini, \emph{A finite group attached to the {L}aplacian of a
  graph}, Discrete Math. \textbf{91} (1991), no.~3, 277--282. \MR{1129991
  (93a:05091)}

\bibitem[LP10]{LevinePropp}
Lionel Levine and James Propp, \emph{What is {$\dots$} a sandpile?}, Notices
  Amer. Math. Soc. \textbf{57} (2010), no.~8, 976--979. \MR{2667495}

\bibitem[Luo11]{Luo}
Ye~Luo, \emph{Rank-determining sets of metric graphs}, J. Combin. Theory Ser. A
  \textbf{118} (2011), no.~6, 1775--1793. \MR{2793609 (2012d:05122)}

\bibitem[Lyo09]{Lyons}
Russell Lyons, \emph{Random complexes and {$l\sp 2$}-{B}etti numbers}, J.
  Topol. Anal. \textbf{1} (2009), no.~2, 153--175. \MR{2541759 (2010k:05130)}

\bibitem[Max09]{Maxwell}
Molly Maxwell, \emph{Enumerating bases of self-dual matroids}, J. Combin.
  Theory Ser. A \textbf{116} (2009), no.~2, 351--378. \MR{2475022
  (2010a:05048)}

\bibitem[May06]{Mayhew}
Dillon Mayhew, \emph{Equitable matroids}, Electron. J. Combin. \textbf{13}
  (2006), no.~1, Research Paper 41, 8 pp. (electronic). \MR{2223516
  (2007c:05047)}

\bibitem[Mer94]{Merris}
Russell Merris, \emph{Degree maximal graphs are {L}aplacian integral}, Linear
  Algebra Appl. \textbf{199} (1994), 381--389. \MR{1274427 (95e:05083)}

\bibitem[Mer05]{Merino-chip-firing}
Criel Merino, \emph{The chip-firing game}, Discrete Math. \textbf{302} (2005),
  no.~1-3, 188--210. \MR{2179643 (2007c:91036)}

\bibitem[Mil66]{Milnor}
J.~Milnor, \emph{Whitehead torsion}, Bull. Amer. Math. Soc. \textbf{72} (1966),
  358--426. \MR{0196736 (33 \#4922)}

\bibitem[MMRW15]{MMRW}
Jeremy~L. Martin, Molly Maxwell, Victor Reiner, and Scott~O. Wilson,
  \emph{Pseudodeterminants and perfect square spanning tree counts}, J. Comb.
  \textbf{6} (2015), no.~3, 295--325.

\bibitem[Moo70]{Moon}
J.~W. Moon, \emph{Counting labelled trees}, From lectures delivered to the
  {Twelfth} {Biennial} {Seminar} of the {Canadian} {Mathematical} {Congress}
  ({Vancouver}, vol. 1969, Canadian Mathematical Congress, Montreal, Que.,
  1970. \MR{0274333 (43 \#98)}

\bibitem[MP95]{MahadevPeled}
N.~V.~R. Mahadev and U.~N. Peled, \emph{Threshold graphs and related topics},
  Annals of Discrete Mathematics, vol.~56, North-Holland Publishing Co.,
  Amsterdam, 1995. \MR{1417258 (97h:05001)}

\bibitem[MR03]{JLM-Factorization}
Jeremy~L. Martin and Victor Reiner, \emph{Factorization of some weighted
  spanning tree enumerators}, J. Combin. Theory Ser. A \textbf{104} (2003),
  no.~2, 287--300. \MR{2019276 (2004i:05070)}

\bibitem[MR05]{JLM-Cyclotomic}
\bysame, \emph{Cyclotomic and simplicial matroids}, Israel J. Math.
  \textbf{150} (2005), 229--240. \MR{2255809 (2007g:05040)}

\bibitem[MR14]{Musiker-Reiner}
Gregg Musiker and Victor Reiner, \emph{The cyclotomic polynomial
  topologically}, J. Reine Angew. Math. \textbf{687} (2014), 113--132.
  \MR{3176609}

\bibitem[MV02]{MV}
Gregor Masbaum and Arkady Vaintrob, \emph{A new matrix-tree theorem}, Int.
  Math. Res. Not. (2002), no.~27, 1397--1426. \MR{1908476 (2003a:05107)}

\bibitem[Oxl11]{Oxley}
James Oxley, \emph{Matroid theory}, second ed., Oxford Graduate Texts in
  Mathematics, vol.~21, Oxford University Press, Oxford, 2011. \MR{2849819
  (2012k:05002)}

\bibitem[Pet09]{Petersson}
Anna Petersson, \emph{Enumeration of spanning trees in simplicial complexes},
  Master's thesis, Uppsala Universitet, 2009.

\bibitem[PPW13]{Primer}
David Perkinson, Jacob Perlman, and John Wilmes, \emph{Primer for the algebraic
  geometry of sandpiles}, Tropical and non-{A}rchimedean geometry, Contemp.
  Math., vol. 605, Amer. Math. Soc., Providence, RI, 2013, pp.~211--256.
  \MR{3204273}

\bibitem[Pr{\"u}18]{Prufer}
H.~Pr{\"u}fer, \emph{Neuer {B}eweis eines {S}atzes {\"u}ber {P}ermutationen},
  Arch. Math. Phys. \textbf{27} (1918), 142--144.

\bibitem[PSFTY09]{WhatIs}
Mohammad~R. Pournaki, Seyed~A. Seyed~Fakhari, Massoud Tousi, and Siamak
  Yassemi, \emph{What is {$\dots$} {S}tanley depth?}, Notices Amer. Math. Soc.
  \textbf{56} (2009), no.~9, 1106--1108. \MR{2568497 (2010k:05346)}

\bibitem[Rei76]{Reisner}
Gerald~Allen Reisner, \emph{Cohen-{M}acaulay quotients of polynomial rings},
  Advances in Math. \textbf{21} (1976), no.~1, 30--49. \MR{0407036 (53
  \#10819)}

\bibitem[RW02]{RW}
Jeffery~B. Remmel and S.~Gill Williamson, \emph{Spanning trees and function
  classes}, Electron. J. Combin. \textbf{9} (2002), no.~1, Research Paper 34,
  24 pp. (electronic). \MR{1928786 (2003g:05067)}

\bibitem[Sco62]{Scoins}
H.~I. Scoins, \emph{The number of trees with nodes of alternate parity}, Proc.
  Cambridge Philos. Soc. \textbf{58} (1962), 12--16. \MR{0136554 (25 \#24)}

\bibitem[SH59]{Simmonard}
M.A. Simmonard and G.F. Hadley, \emph{The maximum number of iterations in the
  transportation problem}, Naval Research Logistics Quarterly \textbf{6}
  (1959), 125--129.

\bibitem[Sta82]{Stanley-LDE}
Richard~P. Stanley, \emph{Linear {D}iophantine equations and local cohomology},
  Invent. Math. \textbf{68} (1982), no.~2, 175--193. \MR{666158 (83m:10017)}

\bibitem[Sta93]{Stanley-Decomp}
\bysame, \emph{A combinatorial decomposition of acyclic simplicial complexes},
  Discrete Math. \textbf{120} (1993), no.~1-3, 175--182. \MR{1235904
  (94k:55027)}

\bibitem[Sta96]{Stanley-CCA}
\bysame, \emph{Combinatorics and commutative algebra}, second ed., Progress in
  Mathematics, vol.~41, Birkh\"auser Boston, Inc., Boston, MA, 1996.
  \MR{1453579 (98h:05001)}

\bibitem[Sta99]{EC2}
\bysame, \emph{Enumerative combinatorics. {V}ol. 2}, Cambridge Studies in
  Advanced Mathematics, vol.~62, Cambridge University Press, Cambridge, 1999,
  With a foreword by Gian-Carlo Rota and appendix 1 by Sergey Fomin.
  \MR{1676282 (2000k:05026)}

\bibitem[Sta07]{Stanley-HA}
\bysame, \emph{An introduction to hyperplane arrangements}, Geometric
  combinatorics, IAS/Park City Math. Ser., vol.~13, Amer. Math. Soc.,
  Providence, RI, 2007, pp.~389--496. \MR{2383131}

\bibitem[SW10]{SuWagner}
Yi~Su and David~G. Wagner, \emph{The lattice of integer flows of a regular
  matroid}, J. Combin. Theory Ser. B \textbf{100} (2010), no.~6, 691--703.
  \MR{2718687 (2012a:05066)}

\bibitem[Syl57]{Sylvester}
J.J. Sylvester, \emph{On the change of systems of independent variables},
  Quart. J. Math. \textbf{1} (1857), 42--56, Collected Mathematical Papers,
  vol. 2, Cambridge, 1908, pp. 65-85.

\bibitem[Tut79]{Tutte-self-dual}
W.~T. Tutte, \emph{On the spanning trees of self-dual maps}, Second
  {I}nternational {C}onference on {C}ombinatorial {M}athematics ({N}ew {Y}ork,
  1978), Ann. New York Acad. Sci., vol. 319, New York Acad. Sci., New York,
  1979, pp.~540--548. \MR{556066 (81c:05031)}

\bibitem[Wil96]{Wilson}
David~Bruce Wilson, \emph{Generating random spanning trees more quickly than
  the cover time}, Proceedings of the {T}wenty-eighth {A}nnual {ACM}
  {S}ymposium on the {T}heory of {C}omputing ({P}hiladelphia, {PA}, 1996), ACM,
  New York, 1996, pp.~296--303. \MR{1427525}

\end{thebibliography}
\end{document}